\documentclass[12pt]{article}
\usepackage[margin=1in]{geometry}  % set the margins to 1in on all sides
\usepackage{graphicx}              % to include figures
\usepackage{amsmath}               % great math stuff
\usepackage{multirow}
\usepackage{amsfonts}              % for blackboard bold, etc
\usepackage{amsthm,float}                % better theorem environments
\usepackage{parskip}
\usepackage[square,sort,comma,numbers]{natbib}
\usepackage[doublespacing]{setspace}

\usepackage{tocloft}
\newcounter{dummy} 
\numberwithin{dummy}{section}
\theoremstyle{definition}
\newtheorem{example}[dummy]{Example}
\newtheorem{definition}[dummy]{Definition}
\newtheorem{theorem}[dummy]{Theorem}

\newtheorem{lemma}[dummy]{Lemma}
\newtheorem{remark}[dummy]{Remark}\newtheorem*{pf}{Proof}

\begin{document}
\nocite{*}

\title{Connected Components of Affine Primitive Permutation Groups}
\author{Haval M. Mohammed Salih} 

\maketitle

\begin{abstract}
For a finite group  $G$, the Hurwitz space $\mathcal{H}^{in}_{r,g}(G)$ is the space of genus $g$ covers of the Riemann sphere with $r$ branch points and the monodromy group $G$.

In this paper, we give a complete list of primitive genus one systems of affine type. That is, we assume that $G$ is a primitive group of affine type. Under this assumption we determine the braid orbits on the suitable Nielsen classes, which is equivalent to finding connected components in $\mathcal{H}^{in}_{r,1}(G)$.
Furthermore, we give a new algorithm for computing large braid orbits on Nielsen classes. This algorithm utilizes a correspondence between the components of $\mathcal{H}^{in}_{r,1}(G)$ and $\mathcal{H}^{in}_{r,1}(M)$, where $M$ is the point stabilizer in $G$.
\end{abstract}

\section{Introduction}
Suppose that $X$ is a compact connected Riemann surface of genus $g$ and that $\mu \colon X\longrightarrow \mathbb{P}^1$ is an indecomposable meromorphic function where $\mathbb{P}^1=\mathbb{C}\cup \{\infty\}$ is the Riemann sphere. For every meromorphic function, there is a number $n$ such that the fiber $\mu^{-1}(p)$ is of size $n$ for all but finitely many points $p\in \mathbb{P}^1$. The number $n$ is called the degree of $\mu$. The points $p$ where $\mu^{-1}(p)<n$ are called branch points of $\mu$. Let $B=\{b_1,...,b_r\}\leq \mathbb{P}^1$ be a finite subset of branch points of $\mu$. Label the points in $\mu^{-1}(p)$ by $\{x_1,...,x_n\}$. The function $\mu$ is not a covering because of ramification points. Thus the restriction function $\mu\colon X\setminus {\mu^{-1}(B)}\rightarrow \mathbb{P}^1\setminus {B}$ is a covering \cite{MR1326604}. The fundamental group $\pi_1(\mathbb{P}^1 \setminus{B},p)$ acts on $\mu^{-1}(p)$ by path lifting \cite{MR1405612}. It induces a group homomorphism $\rho \colon \pi_1(\mathbb{P}^1 \setminus{B},p)\rightarrow S_n$. The image $Mon(X,\mu)$ of $\rho$ is called the monodromy group of $\mu$. It is unique up to conjugacy in $S_n$. The monodromy group $Mon(X,\mu)$ is the Galois group associated to the Galois closure of the extension $C(X)/C(\mathbb{P}^1)$. Since $X$ is connected, then $Mon(X,\mu)$ is a transitive subgroup of $S_n$. Furthermore, $\pi_1(\mathbb{P}^1 \setminus{B},p)$ is generated by all homotopy classes of loops $\gamma_i$ winding once around the point $b_i$ for all $i$. The only relation satisfied by the $\gamma_i$ is $\gamma_1\cdot \gamma_2\cdot...\cdot \gamma_r=1$. Applying $\rho$ to the canonical generators of $\pi_1(\mathbb{P}^1 \setminus{B},p)$ gives the generators of a product one generating tuple in $G$. Simplifying notation we set $x_i=\rho(\gamma_i)$, $1\leq i\leq r$ and $G=Mon(X,\mu)$. The following are true:

\begin{equation} G=\langle x_1,x_2,...,x_r\rangle 
\label{rh1}
\end{equation} 

 \begin{equation} \prod_{i=1}^r {x_i}=1 , \ x_i\in G^{\#}=G\setminus\{1\},\ i=1,...,r. \label{rh2}\end{equation} 
 
\begin{equation}
\sum_{i=1}^r ind\, x_i=2(n+g-1)
\label{rh3}
\end{equation} 
where $ind \,x_i$ is the minimal number of 2-cycles needed to express $x_i$ as a product. Equation (\ref{rh3}) is known as the Riemann-Hurwitz formula. It gives a numerical relation among $g, n$ and $r$. If $C_i$ denotes the conjugacy class of $x_i$. Then the multi set of non trivial conjugacy classes $\bar{C}=\{C_1,...,C_r\}$ in $G$ is called the ramification type of the cover $\mu$.

In light of the above we say that a transitive subgroup $G\leq S_n$ is a genus $g$ group if there exist $x_1,...,x_r\in G$ satisfying (\ref{rh1}), (\ref{rh2}) and (\ref{rh3}) and we call $(x_1,...,x_r)$ the genus $g$ system of $G$. If the action of $G$ on $\{1,...,n\}$ is primitive, we call $G$ a primitive genus $g$ group and $(x_1,...,x_r)$ a primitive genus $g$ system. 

The question is what are possible groups $G$ can occur. For instance, the group $G=\mathbb{Z}_p$ appears while taking the map $\mu_p:\mathbb{P}^1\rightarrow \mathbb{P}^1$ defined by $\mu_{p,a}(z)=(z-a)^p$ is ramified at $a$ and $\infty$.

Our study relates to a conjecture made by Guralnick and Thompson in 1990, in \cite{MR1055011}. In this paper they conjectured that the set $\mathcal{E}^{\ast}(g)$ of possible isomorphism classes of composition factors of $G$, which are neither cyclic nor alternating, is finite for all $g\geq 0$ \cite{MR1055011}. In 2001 the conjecture was proved by Frohardt and Magaard \cite{MR1935417}. The proof of the conjecture shows that we can determine $\mathcal{E}^{\ast}(g)$ explicitly for $g\leq 2$. As the conjecture is now a theorem, these sets are finite.

By assumption $\mu$ is indecomposable in the sense that $\mu$ can not be factorized non-trivially as $\mu=\mu_1\circ \mu_2$, where $\mu_1$ and $\mu_2$ are non-constant functions. So the fact that $\mu$ is indecomposable implies that the monodromy group $G$ acts primitively on the fiber of generic point. In this case the structure of $G$ is explicitly organized around the Aschbacher and O'Nan-Scott Theorem \cite{MR772471}.

\begin{theorem}\label{AC}\cite{MR772471}
Suppose that $G$ is a finite group and $M$ is a maximal subgroup of $G$ such that $$\bigcap_{g\in G} {M^g}=1$$
Let $S$  be a minimal normal subgroup of $G$, let $L$ be a minimal normal subgroup of $S$, and let $\Delta=\left\{L=L_1,L_2,...,L_m\right\}$
be the set of the $G$-conjugates of $L$. Then $L$ is simple, $S=\langle L_1,...,L_r\rangle$, $G=MS$ and furthermore  either
\begin{description}
\item [(A)] $L$ is of prime order $p$; \\
or $L$ is non abelian simple group and one of the following hold:
\item [(B)] $F^{\ast}(G)=S\times R$, where $S\cong R$ and $M\cap S=1$;
\item [(C1)] $F^{\ast}(G)=S$ and $M\cap S=1$;
\item [(C2)] $F^{\ast}(G)=S$ and $M\cap S\neq 1=M\cap L$;
\item [(C3)] $F^{\ast}(G)=S$ and $M\cap S=M_1\times M_2\times\cdots M_t$, where $M_i=M\cap L_i$, $1\leq i\leq t$.
\end{description}
\end{theorem}

In cases (B) and (C1), Shih \cite{MR1129542}, and Guralnick and Thompson \cite{MR1055011} respectively, showed that there are no primitive genus 0 systems. In case (C2) Aschbacher \cite{mr1}, showed that $F^{\ast}(G)=A_5\times A_5$ in all genus 0 examples. In case (C3), $L_i$ is of Lie type of rank one all elements of $\mathcal{E}^{\ast}(0)$ and $\mathcal{E}^{\ast}(1)$ were determined by Frohardt, Guralnick and Magaard \cite{MR1935417}. Furthermore they showed $t\leq 2$. In \cite{MR197} they showed that if $t=1$, $L_i$ is classical and $L_i/M_i$ is a point action, then $n=[L_i,M_i]\leq 10,000$. That result together with results of Aschbacher, Guranlick and Magaard \cite{MR196} show that if $t=1$ and $L_i$ is classical, then $[L_i,M_i]\leq 10,000$ .

In the case (A), where $F^{\ast}(G)$ is abelian and which we refer to as the affine case, Guralnick and Thompson \cite{MR1055011}, showed that there are only finitely many simple groups occurring as composition factors of a primitive genus $0$ group. Furthermore, Neubauer \cite{MR2716273}, showed that there are only finitely many simple groups occurring as composition factors of a primitive genus 0 and 1 group. Finally, Magaard, Shpectorov, and Wang \cite{MR2953829}, produced a complete list of affine primitive genus 0 groups. This case was first considered by Neubauer in his PhD thesis for a genus 1 system. He classified primitive genus 1 systems up to signature. We are interested in this case. Our goal is to classify primitive genus 1 systems up to diagonal conjugation and braiding.

The equivalence classes of $G$-covers $X$ of $\mathbb{P}^1$ with $r$ branched points are called a Hurwitz space and denoted by $\mathcal{H}^{in}_{r,1}(G)$ where $in$ denotes an inner automorphism of $G$. Note that $X$ is a Riemann surface of genus 1 and for the rest $\mathcal{H}^{in}_r(G)$ denotes $\mathcal{H}^{in}_{r,1}(G)$ unless otherwise stated.

Hurwitz shows that the Hurwitz space of simple covers is connected. Also he showed that the connectedness of the Hurwitz space by considering every curve admits a simple cover of $\mathbb{P}^1$.

The Riemann Existence Theorem tells us there is a one to one correspondence between the equivalence classes of product one generating tuples $(x_1,...,x_r)$ of $G$ and the equivalence classes of $G$-covers of type $\bar{C}$ such that $x_i\in C_i$ for $i=1,...,r$.

\begin{theorem}\label{RI}\cite{MR27}
Let $G$ be a finite group and $\bar{C}=\{C_{1},...,C_{r}\}$ be a ramification type. Then there exists a $G$-cover of type $\bar{C}$ if and only if there exists a generating tuple $(x_1,...,x_r)$ of $G$ with $\prod_{i=1}^r x_i=1$ and $x_i\in C_{i}$, for $i=1,...,r$.
\end{theorem}

For any $r$-tuple $(x_1,...,x_r)$ gives a ramification type $\bar{C}$ with $x_i\in C_i$ for $i=1,...,r$. Let $\bar{C}$ be a fixed ramification type, then the subset $\mathcal{H}_r^{in}(G,\bar{C})$ of $\mathcal{H}_r^{in}(G)$ consists of all 
$[P,\phi]$ with admissible surjective map $\phi \colon \pi(\mathbb{P}^1\setminus P,p)\rightarrow G$ sends the conjugacy class $\sum_{p_i}$ to the conjugacy class $C_i$ for $i=1,...,r$.
It is a union of connected components in $\mathcal{H}_r^{in}(G)$.

Here, we study the Hurwitz space $\mathcal{H}_r^{in}(G)$. In particular we focus on the subset $\mathcal{H}_r^{in}(G,\bar{C})$ of $\mathcal{H}_r^{in}(G)$. We try to find the connected components of the Hurwitz space $\mathcal{H}_r^{in}(G)$. To do this, one needs to find corresponding braid orbits.

We present the main theorem of Neubauer's PhD thesis \cite{MR2716273}. It serves as the starting point for our work. 
\begin{theorem}
	\label{A2}
	If $G$ is a primitive genus 1 group of affine type, then
	one of the following holds:
	\begin{enumerate}
		\item $G^{''}=1\text{ and } e\leq 4$.
		\item $p=2 \text{ and } 2\leq e\leq 8$.
		\item $p=3\text{ and } 2\leq e\leq 4$.
		\item $p=5 \text{ or } 7 \text{ and } e\leq 3$.
		\item $p=11 \text{ and } e=2$.
	\end{enumerate}
\end{theorem}

The groups which satisfy 1. in Theorem \ref{A2} are well known \cite{MR2716273} and called Frobenius groups. The other cases in Theorem \ref{A2} will studied.
Our main result, Theorem \ref{r1} gives the complete classification of primitive genus one systems of affine type when $G^{''}\neq 1$.

Appendix $A$ contains tables representing the results of our computation of primitive genus one systems in affine groups satisfying Theorem \ref{A2}, (2)-(5).   

\begin{theorem}\label{r1} 
	Up to isomorphism, there exist exactly 85 affine primitive genus one groups that satisfy Theorem \ref{A2}, (2)-(5). The corresponding primitive genus one groups are enumerated in Tables A.2-17.
\end{theorem}  

Throughout this paper, we will assume that $G$ is a primitive, faithful permutation group on a finite set $\Omega$, $|\Omega|=n$ and $G$ has an abelian minimal normal subgroup $V$. In other words $G$ is a group satisfying (A) of Theorem \ref{AC}. Choose $\omega\in \Omega$ and let $M=G_{\omega}$. As $G$ is primitive, $M$ is a maximal subgroup of $G$. Then 
\begin{equation*}
G=VM, M\cap V=1,  V=C_G(V), |V|=n=p^e \tag{*}
\end{equation*}
for some prime $p$ and for a natural number $e$.

In this paper, we introduce a new algorithm, called the projection-fiber algorithm to compute braid orbits of big length on Nielsen classes and we prove some results related to it. This method applies where $G$ is an affine group. Our idea is to send the tuple $t$ in $G$ via the canonical group homomorphism $\pi\colon G\rightarrow M$ to the tuple $\bar{t}$ in the quotient $M$. This idea is quite useful because the size of the structure constant becomes smaller. An application of this algorithm is the classification of the primitive genus one systems of affine type. That is we find the connected components $\mathcal{H}_r(G,\bar{C})$ of $G$-curves $X$ such that $g(X/G)=0$. In our situation, the computation shows that there is exactly 6000 braid orbits of primitive genus 1 systems of affine type with $G^{''}\neq 1$. The degree and the number of the branch points are given in Table \ref{t1}. This completes the work of Neubauer on the affine case of the Guralnick-Thompson conjecture. Another consequence of this algorithm is that it often gives a one to one correspondence between the number of generating braid orbits of $G$ and $M$ for some groups $G$. That is, the number of components of $\mathcal{H}^{in}_r(G)$ is equal to the number of components of $\mathcal{H}^{in}_r(M)$ for some $G$. For instance this is true for $G=ASL(5,2)$ when $r\geq 4$.

\begin{table}[ht]
	\caption{ Affine Primitive Genus One Systems: Number of Components } % title of Table
	\centering % used for centering table
	\scriptsize\begin{tabular}{|c | c | c| c| c| c| c|c | c|c|c|} % centered columns (4 columns)
		\hline\hline %inserts double horizontal lines
		$Degree$  &  \shortstack{ {\#}Group Iso \\ types}    & {\#} RTs   &  \shortstack{ {\#} comp's\\ $r=3$}   &   \shortstack{ {\#} comp's\\ $r=4$}  &  \shortstack{ {\#} comp's\\ $r=5$}   &   \shortstack{ {\#} comp's\\ $r=6$}    &   \shortstack{ {\#} comp's\\ $r=7$}  &  \shortstack{ {\#} comp's\\ $r=8$}  & \shortstack{ {\#} comp's\\ total}  \\
		\hline % inserts single horizontal line
		128  & 1   & 2 & 2 & -  & -   & -   &  -    &  - & 2\\
		\hline % inserts single horizontal line
		64 &  24  &  114  &  738  & 19 &  -  &  -  & -  &  - &757 \\
		\hline % inserts single horizontal line
		32 &  1  &  131 &  2247 &  30 & 3  & - & - & - & 2280\\
		\hline % inserts single horizontal line
		16 &  18  &  599 &  2015 &  273  & 94  & 18 & 1 & - &2409 \\
		\hline % inserts single horizontal line
		8  &  2  &  134 &  64 &  71  & 28  & 13 & 4 & 1 &  181 \\
		
		\hline % inserts single horizontal line
		9  &  5  &  49 &  26 &  25  & 9  & 3 & - & - & 63\\
		\hline % inserts single horizontal line
		27  &  7  &  53 &  119 &  14  & 2  & - & - & - & 135  \\
		\hline % inserts single horizontal line
		81  &  14  &  37 &  71 &  4  & -  & - & - & - & 75\\
		
		\hline % inserts single horizontal line
		25  &  5  &  17 &  34 &  1  & -  & - & - & - & 35\\
		\hline % inserts single horizontal line
		125  &  2  &  5 &  24 &  -  & -  & - & - & - & 24\\
		\hline % inserts single horizontal line
		49  &  4  &  15 &  14 &  1  & -  & - & - & - & 15\\
		\hline % inserts single horizontal line
		343  &  1  &  12 &  12 &  -  & -  & - & - & - & 12\\
		\hline % inserts single horizontal line
		121  &  1  &  4 &  20 &  -  & -  & - & - & - & 20\\
		\hline % inserts single horizontal line
		Totals    &  85  &  1172 &  5386 &  438  & 136  & 34 & 5 & 1 & 6000\\
		\hline % inserts single horizontal line
	\end{tabular}
	\label{t1}
\end{table}

This paper is organized as follows. Section 2 is a background results and here we explain the relationship between connected components in Hurwitz spaces and braid orbits on Nielsen classes. Section 3 is devoted to order and label conjugacy classes of $G$ according to a certain rules to distinguish classes represented by elements of equal order. In Section 4 the projection-fiber algorithm is introduced and several results are given about it. Finally, we will give some examples to explain our algorithm. In Section 5 the methodology are give which we use to compute braid orbits for primitive genus 1 systems of affine type.

\section{Preliminaries}
We begin this section with a formal definition of the Artin braid group.
\begin{definition}
For $r\geq 2$, the Artin braid group $B_r$ is generated by $r-1$ elements $\sigma_1,\sigma_2,...,\sigma_{r-1}$ that satisfy the following relations:
\begin{equation}\sigma_i\sigma_j=\sigma_j\sigma_i 
\end{equation} for all $i,j=1,2,...,r-1$ with $|i-j|\geq 2$, and \begin{equation}\sigma_i\sigma_{i+1}\sigma_i=\sigma_{i+1}\sigma_{i}\sigma_{i+1}
\end{equation} for $i=1,2,...,r-2$. These relations are known as the braid relations.
\end{definition}

The braid $\sigma_i$ acts on generating tuples $x=(x_1,...,x_r)$ of a finite group $G$ with $\prod_{i=1}^r{x_i}=1$ as follows:
\begin{equation}
(x_1,...,x_i,x_{i+1},...,x_r)\sigma_i=(x_1,...,x_{i+1},x_{i+1}^{-1}x_ix_{i+1},...,x_r) 
\label{op}
\end{equation} 
for $i=1,...,r-1$.
The braid orbit of $x$ is the smallest set of tuples which contains $x$ and is closed under the operations (\ref{op}).

Applying $\phi \colon \pi(\mathbb{P}^1\setminus P,p)\rightarrow G$ to the canonical generators of $\pi_1(\mathbb{P}^1\setminus P,p)$ gives the generators of a product one generating tuple in $G$ that is, $\phi(\lambda_i)=x_i$. We define $\epsilon_r(G)=\{(x_1,...,x_r): G=\langle x_1,...,x_r\rangle, \prod_{i=1}^r{x_i}=1, x_i\in G^{\#}, i=1,...,r\}$. Let $A\leq Aut(G)$. Then the subgroup $A$ acts on $\epsilon_r(G)$ via sending $(x_1,...,x_r)$ to $(a(x_1),...,a(x_r))$, for $a \in\mathrm{A}$, which is known as the diagonal conjugation. This action commutes with the operations (\ref{op}). Thus $A$ permutes the braid orbits. If $A=Inn(G)$, then it leaves each braid orbit invariant \cite{MR1405612}. Let $\epsilon_r^{in}(G)=\epsilon_r(G)/Inn(G)$.

We now introduce the Nielsen classes in $G$, as follows. For a ramification type $\bar{C}$,\\
$\mathcal{N}(\bar{C})=\{(x_1,...,x_r): G=\langle x_1,...,x_r\rangle, \prod_{i=1}^r{x_i}=1, \exists\sigma\in S_n \text{ such that } x_i\in C_{i\sigma} \text{ for all } i\}.$

Assume $(*)$ holds and since $\phi \colon \pi_1(\mathbb{P}^1\setminus P,p)\rightarrow G$ is a surjective homomorphism and $\pi \colon G\rightarrow M$ is the canonical group homomorphism, then $\eta:=\phi \circ \pi \colon \pi_1(\mathbb{P}^1\setminus P,p)\rightarrow M$ is also a surjective homomorphism which sends the canonical generators of $\pi_1(\mathbb{P}^1\setminus P,p)$ to generators of $M$, say $\eta(\lambda_i)=m_i$. Similarly we can define $\epsilon_r(M) $ and $\mathcal{N}(\bar{\bar{C}})$ in $M$ where $\bar{\bar{C}}=\{\bar{C_1},...,\bar{C_l}\}$ is a ramification type of $m=(m_1,...,m_l)$. It is clear that the size of $\mathcal{N}(\bar{\bar{C}})$ is less than the size of $\mathcal{N}(\bar{C})$.

The topology on $\mathcal{H}_r^{\mathrm{A}}(G)$ is well defined. Let $O_r$ be the set of all $r$-tuples of distinct elements in $\mathbb{P}^1$, equipped with the product topology \cite{MR1119950}.

The next result is account to begin our investigation.
\begin{lemma}\cite{MR1405612}
\label{t11}
The map $\Psi_{\mathrm{A}}:\mathcal{H}_r^{\mathrm{A}}(G) \longrightarrow O_r$, $\Psi_{\mathrm{A}}([P,\phi]))=P$ is covering.
\end{lemma}

The fundamental group $\pi_1(O_r,P_0)=B_r$ acts on $\Psi_{\mathrm{A}}^{-1}(P_0)$ where $P_0=\{1,...,r\}$ is the base point in $O_r$ via path lifting where the fiber $\Psi_{\mathrm{A}}^{-1}(P_0)=\{[P_0,\phi]_{\mathrm{A}}:\phi\colon \pi_1(\mathbb{P}^1\setminus{P_0},\infty)\longrightarrow G \text{ is admissible }\}$. This $\phi$ gives a product one generating tuple $(x_1,...,x_r)$ of $G$. 

\begin{lemma}\cite{MR1405612}
\label{t2}
We obtain a bijection $\Psi_{\mathrm{A}}^{-1}(P_0)\longrightarrow \epsilon_r^{\mathrm{A}}(G)$ by sending $[P_0,\phi]_{\mathrm{A}}$ to the generators $(x_1,...,x_r)$ where $x_i=\phi([\gamma_i])$ for $i=1,...,r$.
\end{lemma} 

The image $\mathcal{N}^{A}(\bar{C})$ of $\mathcal{N}(\bar{C})$ in $\epsilon_r^{\mathrm{A}}(G)$ is the union of braid orbits. 
If $\Psi_{\mathrm{A}}$ in Lemma \ref{t11} restricts to a connected component $\mathcal{H}$ of $\mathcal{H}_r^{\mathrm{A}}(G)$, then Lemma \ref{t2} implies that the fiber in $\mathcal{H}$ over $P_0$ corresponds to the set $\mathcal{N}^{A}(\bar{C})$. This yields a one to one correspondence between connected components of $\mathcal{H}_r^{\mathrm{A}}(\bar{C})$ and the braid orbits on $\mathcal{N}^{A}(\bar{C})$. For $A=Inn(G)$, we see that the above gives a one to one correspondence between connected components of $\mathcal{H}_r^{in}(\bar{C})$ and the braid orbits on $\mathcal{N}(\bar{C})$ \cite{MR1405612}.

\begin{definition}\cite{MR27}
Two generating tuples are braid equivalent if they lie in the same orbit under the group generated by the braid action and diagonal conjugation by $Inn(G)$.
\end{definition} 

That is if two generating tuples lie in the same braid orbit under either the diagonal conjugation or the braid action, then the corresponding covers are equivalent by Riemann's Existence Theorem.

\begin{definition}
Two coverings $\mu_1\colon X_1 \rightarrow \mathbb{P}^1$ and $\mu_2\colon X_2 \rightarrow \mathbb{P}^1$ are equivalent if there exists a homeomorphism $\alpha \colon X_1\rightarrow X_2$ with $\mu_2\alpha=\mu_1$.
\end{definition}
\begin{theorem}\cite{MR1405612}
Two generating tuples are braid equivalent if and only if their corresponding covers are equivalent.
\end{theorem}

To answer whether or not $\mathcal{H}_r(G,\bar{C})$ is connected is still an open problem, both computationally and theoretically. The computation becomes difficult when the length of Nielsen classes grows rapidly. The MAPCLASS package of James, Magaard, Shpectorov and Volklein, is designed to perform braid orbit computations for a given finite group and given type. Few results were known about it. For instance, Clebsch \cite{MR1509816} shows that if $G=S_n$ and let $\textbf{C}=(C,...,C)$ be $r$-tuple consisting of $r$ copies the class $C$ of transpositions, then the corresponding Hurwitz space $\mathcal{H}^{in}_r(G,\textbf{C})$ is connected. Liu and Osserman \cite{MR2464030} generalized this result as follows. If $G=S_n$ and $C_i$ represented by $x_i$ where $x_i$ is a single cycle of length $|x_i|$, then $\mathcal{H}_r(G,\bar{C})$ is connected. Furthermore, Fried \cite{fried2006alternating} shows that if $G=A_n$, $g>0$, and all $C_i$ are represented by 3-cycles  then $\mathcal{H}_r(G,\bar{C})$, has one component if $g_1=g(X/G)=0$ and otherwise it has two components.

The problem of classifying the braid orbits appeared from the study of braid monodromy factorization. The answering is well known for solvable groups. However, there are few results on the classification of braid orbits for non-solvable groups. Ben-Itzhak and Teichen \cite{MR1996396} determine all braid orbits on Nielsen classes for symmetric group of degree $n$. Recently Magaard, Shpectorov, and Wang\cite{MR2953829} determined all braid orbits on Nielsen classes of primitive genus zero systems of affine type.

\section{Ordered and labeling convention of conjugacy classes}
We are going to order and label conjugacy classes according to a certain rule to distinguish classes represented by elements of equal order. This has to do with the fact that, the ordering of the conjugacy classes in a group may not the same for two different runs of GAP. For better consistency we want to establish a canonical order of classes similar to the ATLAS notation. This works as follows.
Let $G$ be an affine primitive permutation group and let $C_1,...,C_n$ be the conjugacy classes of $G$. Assume that $x_i\in C_i$ has order $d_i$ for $i=1,...,n$. The following rules are applied step by step:

\begin{description}
	\item[The Order.] We order the $d_i$, for $i=1,...,n$ in a non-decreasing sequence. 
	\item [The Centralizer size.] If $d_1,...,d_s$ are equal orders for some $s$, then  we compute $|C_G(C_i)|$ for $i=1,...,s$ and ordered it from a decreasing sequence.
	\item [The permutation indices.] If equality holds among some terms in sequence of $|C_G(C_i)|$, then we compute the permutation indices in the natural action (that is, $G$ acts on $p^e$ points) for them and ordered it from a non-decreasing sequence. 
	\item [Power functions.] If equality holds among some terms in sequence of indices, then we take power functions for the representative $x_i$ of conjugacy classes $C_i$ which are correspond to the equally indices, that is $x_i^m$ where $m$ is some order of an element in $G$. Next we are going to check that $x_i^m$ are conjugate to representative elements $y_i$ in different known types in $G$ or not. In Table \ref{s1}, we will see this case for $C_4$ and $C_5$.
\end{description}

According to the above rules, we obtain a sequence. We label the corresponding conguacy class of the first term in sequence by $d_iA$, the corresponding conjugacy class of the second term in sequence by $d_iB$ and so on.

\begin{example}
	In $G=AGL(3,3)$, we look at elements of order six, of which $G$ has eight different conjugacy classes. We denoted by $C_1,C_2,C_3,C_4,C_5,C_6,C_7 \text{ and } C_8$. We first compute the centralizer size and index of each of them and the results appear in Table \ref{s1}. In fact, 6A represents $C_8$ because it has a bigger centralizer size. However, the conjugacy classes $C_3$ and $C_6$ have same centralizer size but the indices are different and so $6B$ must be represent $C_3$, because it is index less than the index of $C_6$ and so on. For the conjugacy classes $C_4$ and $C_5$, we use different technique such as we take the third power of the representative elements in conjugacy classes $C_4$ and $C_5$ which are conjugate to the representative elements in conjugacy classes $2A$ and $2B$ respectively. Therefore $6F$ represent $C_4$ and $6G$ represent $C_5$. 
	
	\begin{table}[ht]
		\caption{Ordering and Labeling Conjugacy Classes} % title of Table
		\centering % used for centering table
		\scriptsize\begin{tabular}{|c | c | c| c|c|c|c|c|c|} % centered columns (4 columns)
			\hline\hline %inserts double horizontal lines
			Conjugacy classes &  $C_1$ & $C_2$ & $C_3$ & $C_4$ & $C_5$ & $C_6$ & $C_7$ & $C_8$  \\
			\hline % inserts single horizontal line
			Size of Centralizer &  36 & 36 & 108 & 18 & 18 & 108 & 18 & 144  \\
			\hline % inserts single horizontal line
			indices &  17 & 18 & 19 & 21 & 21 & 21 & 22 & 22  \\
			\hline % inserts single horizontal line
			After ordering                           \\
			
			\hline % inserts single horizontal line
			ordered &  $C_8$ & $C_3$ & $C_6$ & $C_1$ & $C_2$ & $C_4$ & $C_5$ & $C_7$  \\
			\hline % inserts single horizontal line
			size of centralizer &  144 & 108 & 108 & 36 & 36 & 18 & 18 & 18  \\
			\hline % inserts single horizontal line
			indices &  22 & 19 & 21 & 17 & 18 & 21 & 21 & 22  \\
			\hline % inserts single horizontal line
			power function &  - & - & - & - & - & $C_4^3$ & $C_5^3$ & -  \\
			\hline % inserts single horizontal line
			Types &  6A & 6B & 6C & 6D & 6E & 6F & 6G & 6H  \\
			\hline % inserts single horizontal line
		\end{tabular}
		\label{s1}
	\end{table}
\end{example}

\newpage
\begin{remark}
\begin{enumerate}
\item	Some time it is not useful to take power functions of conjugacy classes because they conjugate to the same element in a group. Also if the conjugacy classes are inverse of each other, then we cannot make distinction among them.
\item 	We label the types as follows. The orders must be a non decreasing and all conjugates elements must be adjacent. For instance if we have the type $(2A,2C,3A,2A)$, then it should be reorder as $(2A,2A,2C,3A)$.
\end{enumerate} 
\end{remark}

\section{Projection-fiber Algorithm}
Assume that (*) holds and let $\pi\colon G\rightarrow M$ be the canonical group homomorphism defined by $\pi(x)=m$ where $x=vm$. Our idea is to send a generating tuple $t$ of $G$ via $\pi$ to the tuple $\bar{t}$ in the quotient $M$. This idea is quite useful because the size of $M$ is less than the size of $G$. In other words, the size of the structure constant becomes smaller.

Let $t=(x_1,...,x_r)$, with $r\geq 4$, be a generating tuple of $G$ of type $\bar{C}$ and $Z(M)=1$. Applying $\pi$ to $t$, we obtain a tuple $\bar{t}=(m_1,...,m_r)$. Note that $\bar{t}$ is generated for $M$. In such situation one can compute the generating braid orbit $O_M=(\tilde{m_1},...,\tilde{m_r})$ for $\bar{t}$ in $M$ by using this function \textbf{GeneratingMCOrbits} in MAPCLASS package. Once this is done, pull the result back by taking pre-images of each element in generating tuple from each braid orbit $O_M$ under the homomorphism $\pi$. In this way, each element from generating tuple from each braid orbits $O_M$ gives the set $L_i=\pi^{-1}(\tilde{m_i})$. The size of each of $L_i$ is equal to the size of $V$. We are interested in these elements in $L_i$ which are conjugate to $x_i$ for $i=1,...,r$. So we can assume that $U_i= L_i\cap x_i$, $x_i\in C_i$ and $U=U_1\times U_2\times ...\times U_r$. Once this is done, we collect these tuples in $U$ which satisfies $\prod_{i=1}^r{u_i}=1$ and $G=\langle u_1,...,u_r\rangle$ in a list which denotes by $\mathcal {G}_l$. This $V$ acts on $\mathcal{G}_l$ via diagonal conjugation. The result of this action is $q$-orbits, where $q\in \mathbb{Z}^+$.

It will have noticed that we are writing the canonical group homomorphsim $\pi$, on the left. That is, we write $\pi(g)$. So we now prove some results related to the algorithm.
\begin{lemma}\label{po1}
Let $\pi:G^r\rightarrow M^r$ be the group homomorphism defined by $\pi(x_1,x_2,...,x_r)=(m_1,...,m_r)$. Then $\pi$ commutes with braid action.
\end{lemma}
\begin{pf}
Let $(x_1,...,x_r)$ be a tuple in $G^r$ and $\sigma_i\in B_r$. Then we have
\begin{align*}
\sigma_i(\pi(x_1,x_2,...,x_r)) &= \sigma_i(m_1,...,m_r)\\
               &=(m_1,...,m_{i+1},m_{i+1}^{-1}m_im_{i+1},...,m_r) \\
               &=(\pi(x_1),...,\pi(x_{i+1}),\pi(x_{i+1}^{-1}x_ix_{i+1}),...,\pi(x_r)) \\
               &=\pi(x_1,...,x_{i+1},x_{i+1}^{-1}x_ix_{i+1},...,x_r) \\
               &=\pi(\sigma_i(x_1,x_2,...,x_r))
               \end{align*}
Hence $\sigma_i\pi=\pi\sigma_i$.
\end{pf}

The next lemma tells us if we have two tuples $t_1$ and $t_2$ which are not braid equivalent in $M$, then $\pi^{-1}(t_1)$ and $\pi^{-1}(t_2)$ are also not braid equivalent in $G$.
\begin{lemma}
Let $G=VM$, where $M$ is a maximal subgroup of $G$ and $V$ is an abelian normal subgroup of $G$. If $\pi:G\rightarrow M$ is the canonical group homomorphism and $(m_1,...,m_r)$ and $(\bar{m_1},...,\bar{m_r})$ are not braid equivalent then $(x_1,...,x_r)$ and $(\bar{x_1},...,\bar{x_r})$ are not braid equivalent in $G$.
\end{lemma}

\begin{pf}
Suppose that $t=(x_1,...,x_r)$ and $\bar{t}=(\bar{x_1},...,\bar{x_r})$. Since $\pi$ is a homomorphism, Lemma \ref{po1} implies $\pi(t^\sigma)=(\pi(t))^\sigma $. So if $\bar{t}=t^\sigma$ for some $\sigma\in B_r$ then
$\pi(\bar{t})=\pi(t^\sigma)=(\pi(t))^\sigma$. Therefore $\pi(\bar{t}) \text{ and } \pi(t)$ are braid equivalent.
\end{pf}

\begin{lemma}
Assume that (*) holds and suppose $m=(m_1,...,m_r)$ is a generating tuple in $M$, and $l(m)=(l(m_1),...,l(m_r))$ is an arbitrary lift of $m$ into $G$. Then either $\langle l(m_1),...,l(m_r)\rangle =G$ or $\langle l(m_1),...,l(m_r)\rangle$ is a complement to $V$.
\end{lemma}

\begin{pf}
It is clear that $G_0=\langle l(m_1),...,l(m_r)\rangle $ is a subgroup of $G$. Let $H=V\cap G_0$, then $H$ is normal in $G_0$ since $V$is normal in $G$. It is also normal in $V$ because $V$ is abelian. Since $G=VG_0$. It follows that $H$ is normal in $G$ and thus either $H=1$ or $H=V$. We are done.
\end{pf}

\begin{lemma}
\begin{enumerate}
\item The number of components of $\mathcal{H}^{in}_{r,0}(M)$ is equal to the number of components of $\mathcal{H}^{in}_{r,0}(G)$, where $r\geq 4$.
\item If $G\neq 2^4.S(16)$, then  $\mathcal{H}^{in}_{r,1}(G,\bar{C})$ is connected where $r\geq 5$. 
\end{enumerate} 
\label{hh1}
\end{lemma}
The proof of Lemma \ref{hh1} 1., can be found in \cite{MR2953829}. Also we apply our algorithm to find the number of components of $\mathcal{H}^{in}_{r,0}(M)$.

\begin{definition}
Let $G$ be a group act on $\Omega$. For $x\in G$ define $Fix\,x=\{\omega\in \Omega : \omega x=\omega\}$ and $f(x)=|Fix\,x|$.
\end{definition}

Recall that $G=VM$ is the semi-direct product of $V$ by $M$ which acts on $n=p^e$ points and $M$ acts on $n-1=p^e-1$ points.
\begin{lemma}
Let $G=VM=\langle x_1,...,x_r\rangle$ with $x_i=v_im_i$ for $i=1,...,r$ and $\pi:G \rightarrow G/V\cong M$ be a canonical group homomorphism define by $\pi(x_i)=m_i$. Let $\bar{x}=(x_1,...,x_r)$ be a generating tuple of $G$ with genus $g$. Then the image $\pi(\bar{x})$ is a generating tuple of $M$ with genus $g^*$ where $g^*\leq g+1$. More precisely, if $f(x_i)>0$, for $i=1,...,r$ then $g^*=g+1$. If $f(x_i)=0$ for some $i\in\{1,...,r\}$ then $g^*=g+1-k$ for some $k\in \mathbb{Z}^+$.
\end{lemma}

\begin{pf}
Since $M=\pi(G)=\pi(\langle x_1,...,x_r\rangle)=\langle \pi(x_1),...,\pi(x_r)\rangle=\langle m_1,...,m_r\rangle$. Also $1=\pi(1_G)=\pi(\prod_{i=1}^r{x_i})=\prod_{i=1}^r{\pi(x_i)}=\prod_{i=1}^r{m_i}$.

In fact, if $f(x_i)>0$ for $i=1,...,r$ then $x_i^x=m_i$ for some $x\in G$ and $i=1,...,r$ which implies that $ind \,x_i=ind \,m_i$. From the Riemann-Hurwitz formula we have that $2(n+g-1)=2(n-1+g^{*}-1)$. Hence $g^{*}=g+1$.

If $f(x_i)=0$ for some $i\in\{1,...,r\}$ then $ind \,x_i\geq ind \,m_i$. It follows that $2(n+g-1)\geq 2(n-1+g^{*}-1)$ there exist a positive even integer number $l$ such that $2(n+g-1)= 2(n-1+g^{*}-1)+l$, we can write $l=2k$ for some $k\in \mathbb{Z}^+$. Hence $g^{*}=g+1-k$.
\end{pf}

\begin{remark}
\label{ha1}
In our situation we conclude that the number of generating braid orbits for a quotient tuple $\bar{t}$ in $M$ is bounded by the number of generating braid orbits for an original tuple $t$ in $G$ which is also bounded by the number of generating braid orbits for a lifting back generating tuple from each $O_M$ to $G$. To be precise, we can write $q=q_1q_2$, where $q_1,q_2\in \mathbb{Z}^+$.  One of them is the number of generating braid orbits in $G$ and the others multiply the length of generating braid orbits in $M$ is the length of braid orbits in $G$.
\end{remark}

From Remark \ref{ha1}, we see that there is two choices either $q_1$-orbit of length $q_2\times length(O_M)$ or $q_2$-orbit of length $q_1 \times length(O_M)$. We can decide which one appears just by running a function \textbf{AllMCOrbits} in MAPCLASS package for a while.

If $q=1$, then there is a one to one correspondence between the number of generating braid orbits in $M$ and the number of generating braid orbits in $G$. That is the number of generating braid orbits in $G$ and $M$ are equal and the length of generating braid orbits in $G$ and $M$ are also equal.

The question comes up if $q\neq 1$, that is the lifting back tuple give at least two tuples. Are they equivalent? The answer it may be equivalent or not equivalent. For more detail we will give the following example.

\begin{example}
In $G=ASL(3,2)$, the type $(3A,3A,3A,3A)$ is generating type. We send tuple which correspondence this type via the canonical group homomorphism $\pi \colon G\longrightarrow M$ to the quotient tuple in $M\cong SL(3,2)$. We can compute the generating braid orbit for quotient tuple in $M$. We obtain two generating braid orbits $O_M(1)$ and $O_M(2)$ of lengths 90 and 144 respectively. After lifting the generating tuple from $O_M(1)$ back to $G$ and $V$ acts diagonal on lifting generating tuples of $G$. We obtain $q=2=1\times 2$. It gives 1 orbit of length $90\times 2$. Similarly for $O_M(2)$. However, in this case we have 2 orbits with length 144. 
\end{example}
We observe that the lifting generating tuple from $O_M(1)$ gives two generating tuples for $G$ which are equivalent. However the lifting generating tuple from $O_M(2)$ gives two generating tuples for $G$ which are not equivalent

Another useful feather of our algorithm concerns the braid orbit computation time in a whole group $G$ as follows:
Let $\pi:G\rightarrow M$ be the canonical group homomorphism and let $t_1,...,t_s$ be generating tuple of $G$ whose $\pi(t_i)=t$. Of course $t$ generates $M$. We compute braid orbits $O_M$ for $t$ in $M$ and then lifting back the generating tuple from $O_M$ to $G$. In this way we achieve the number of braid orbits for $t_1,...,t_s$. We will give the following example.

\begin{example}
In $ASL(5,2)$, and take 2 types $(2D,2D,2E,12B)$ and $(2D,2D,2D,12C)$ whose images of corresponding tuples of these 2 types is a tuple $t$ of type $(2B,2B,2B,12A)$ and it has 1 orbit of length 720. The lifting generating tuple of this orbit gives 1 orbit for each of these 2 types in $G$ and the length of each is 720.
\end{example}

\begin{definition}
A $G$-cover $X$ of Riemann sphere is called $(M,r,g)$ full cover of Riemann sphere if the number of generating braid orbits in $M$ is equal to the number of generating braid orbits in $G$.
\end{definition}

The following example illustrating the application of the projection-fiber algorithm.
\begin{example}
The group $G=2^4.PSL(4,2)$ acting on 16 points of $\mathbb{F}^4_2$ which is semi-direct product of the vector space $V=\mathbb{F}^2_4$ with $PSL(4,2)$, that is $G=V PSL(4,2)$, we have only one generating type of length 7, which corresponds to the ramification type $\bar{C}=(2B,2B,2B,2B,2B,2D,2D)$ and the constant structure for $\bar{C}$ is 1,137,259,549,440. The relevant number is the size of $G$ divides  structure constant which is equal to $\frac{1137259549440}{322560}\approx 3525730$ as it is an estimate for the sum of the lengths of all \textbf{GeneratingMCOrbits}. The computation of this tuple may be impossible directly. Apply the canonical group homomorphism $\pi$ to $\bar{C}$ gives the ramification type $\bar{\bar{C}}=(2A,2A,2A,2A,2A,2B,2B)$ in the quotient, which is isomorphic to $PSL(4,2)$. The structure constant of $\bar{\bar{C}}$ is 29,632,277,430. The relevant number is the size of $PSL(4,2)$ divides structure constant which is equal to $\frac{29632277430}{20160}\approx 1469855$.  Now we compute the generating orbits for corresponding tuple of $\bar{\bar{C}}$ in the quotient $PSL(4,2)$. After we found the generating braid orbit for the quotient tuple in $PSL(4,2)$ and then pull the result back to $G$. We will see that $q=1$. Note that the size of $PSL(4,2)$ divides the total number of tuples (this number appears when we run the \textbf{GeneratingMCOrbits}) , that is $\frac{18192384000}{20160}=902400$.
\end{example}

\begin{table}[H]
\caption{Summarize of Example } % title of Table
\centering % used for centering table
\scriptsize\begin{tabular}{|c | c | c| c| c| c|c|} % centered columns (4 columns)
\hline\hline %inserts double horizontal lines
   Group       &  Type                    &  {\#} of orbits   &  Largest length of orbit  &  Size of constant structure & Time spent  \\
\hline %inserts double horizontal lines
$PSL(4,2)$     & (2A,2A,2A,2A,2A,2B,2B)   &   1 &  902400         &  29,632,277,430    & 1782 minutes\\
\hline % inserts single horizontal line
$2^4.PSL(4,2)$ & (2B,2B,2B,2B,2B,2D,2D)   &  1  &  902400         & 1,137,259,549,440 & $>$ 17280 minutes \\
\hline % inserts single horizontal line
\end{tabular}
\end{table}

\section{Methodology: Listing primitive genus one systems of affine type}
We are presenting our results in Tables A.2-A.17 \cite{MR227}. To obtain these tables we needed to do the following steps:

\begin{itemize}
	\item We extract all primitive permutation group $G$ by using the GAP function \\ AllPrimitiveGroups(DegreeOperation,$p^e$).  Furthermore, we check the order of the socle $V$ of each of those primitive groups to know which them are affine.
	\item For every affine group $G$, compute the conjugacy class representatives and permutation indices on $|V|=p^e$ points.
	\item For given  $p,e,g$ and $G$ we use the GAP function RestrictedPartions to compute all possible ramification types satisfying the Riemann-Hurwitz formula.
	\item For each conjugacy class representative $x$ compute $dim_V(x)$ and use Scott's Theorem to eliminate those types from the previous step which cannot possibly act irreducibly on $V$, that is, they cannot generate a primitive group. 
	\item For each conjugacy class representative $x$ compute $dim_V(x)$ and the number of fixed points of $x$ and use the corollary of Scott's Theorem (Corollary 3.21 in \cite{MR2716273}) to eliminate those types from the previous step which cannot possibly act irreducibly on $V$, that is, again they cannot generate a primitive group. 
	\item Compute the character table of $G$ if possible and remove those types which have zero structure constant.
	\item For each of the remaining types of length greater than or equal to 4, we use MAPCLASS package to compute braid orbits if possible. Otherwise when the length is too big, we use the projection-fiber algorithm. For tuples of length 3 determine braid orbits via double cosets \cite{MR27}. 
	\item  We show that $AGL(8,2)$ possesses no primitive genus 1 systems by using Lemma 3.2 in \cite{MR1865973} and Lemma 3.10 in \cite{MR27}.
\end{itemize}

\section{Acknowledgment}
The results presented in this paper are part of a thesis which the author completed in 2014 at the University of Birmingham under direction of Kay Magaard and Sergey shpectorov and to whom the author wishes to express his thanks for the encourgement and interest in the completion of this project. Especially, it is dedicated to Kay, who died in 2018. 

\appendix
\section{Appendix}
%\begin{singlespace}
%	\include{cover}
%\end{singlespace}
Note that N.O means number of orbits, L.O means largest length of the orbit and GOS means Genus one System.
\begin{table}[H]
	\caption{GOSs for Affine Primitive Groups of Degree 121} % title of Table
	\centering % used for centering table
	\scriptsize\begin{tabular}{|c | c | c| c| c| c| c|} % centered columns (4 columns)
		\hline\hline %inserts double horizontal lines
		group &  ramification type & N.O & L.O &ramification type & N.O & L.O \\
		\hline % inserts single horizontal line
		$ASL(2,11):2$  &(2B,3A,10A) & 5 & 1& (2B,3A,10B) & 5 &1\\
		&(2B,3A,10C) & 5 & 1& (2B,3A,10D) & 5 &1\\
		\hline % inserts single horizontal line
	\end{tabular}
\end{table}

\begin{table}[H]
	\caption{GOSs for Affine Primitive Groups of Degree 49 and 343} % title of Table
	\centering % used for centering table
	\scriptsize\begin{tabular}{|c | c | c| c| c| c| c|} % centered columns (4 columns)
		\hline\hline %inserts double horizontal lines
		group &  ramification type & N.O & L.O &ramification type & N.O & L.O \\
		\hline % inserts single horizontal line
		$7^2:S_3 $    &(2A,3A,14A) & 1 & 1 & (2A,3A,14B) & 1 & 1\\
		&(2A,3A,14C) & 1 & 1 & (2A,3A,14D) & 1 & 1\\
		&(2A,3A,14E) & 1 & 1 & (2A,3A,14F) & 1 & 1 \\
		\hline % inserts single horizontal line
		$7^2:3\times D(2*6)$    & (2A,6C,6H) & 1 & 1 & (2A,6D,6I)  & 1  & 1\\
		&(2B,6C,6G) & 1 & 1 & (2B,6D,6F)  & 1  & 1\\
		\hline % inserts single horizontal line
		$7^2:D(2*6)$          &(2A,2B,2C,3A) & 1 & 24 &  &  &\\
		\hline % inserts single horizontal line
		$7^2:3\times(Q_8:3)$ &(3D,3F,6C) & 1 & 1 &(3C,3F,6E)  &1  & 1\\
		&(3B,3G,6D) & 1 & 1 &(3A,3G,6F)  &1  & 1\\
		\hline % inserts single horizontal line
		$ 7^3:PSL(2,7)$ &(2A,3A,7U) & 1 & 1& (2A,3A,7V) & 1 &1\\
		&(2A,3A,7K) & 1 & 1& (2A,3A,7L) & 1 &1\\
		&(2A,3A,7M) & 1 & 1& (2A,3A,7N) & 1 &1\\
		&(2A,3A,7O) & 1 & 1& (2A,3A,7P) & 1 &1\\
		&(2A,3A,7Q) & 1 & 1& (2A,3A,7R) & 1 &1\\
		&(2A,3A,7S) & 1 & 1& (2A,3A,7T) & 1 &1\\
		\hline % inserts single horizontal line
	\end{tabular}
\end{table}

\begin{table}[H]
	\caption{GOSs for Affine Primitive Groups of Degree 25 and 125} % title of Table
	\centering % used for centering table
	\scriptsize\begin{tabular}{|c | c | c| c| c| c| c|} % centered columns (4 columns)
		\hline\hline %inserts double horizontal lines
		group &  ramification type & N.O & L.O &ramification type & N.O & L.O \\
		\hline % inserts single horizontal line
		AGL(2,5) & (2B,4F,24A) & 1 & 1& (2B,4F,24B) & 1 & 1 \\
		& (2B,4E,24C) & 1 & 1& (2B,4E,24D) & 1 & 1 \\
		& (3A,4E,4F) & 6 & 1 &  &  &  \\
		\hline % inserts single horizontal line
		ASL(2,5) & (3A,4A,5C) & 2 & 1& (3A,4A,5B) & 2 & 1 \\
		\hline % inserts single horizontal line
		$5^2:4\times D(2*3)$   & (2B,4C,12B) & 1 & 1 & (2B,4D,12A) & 1 & 1 \\
		& (2A,4C,12A) & 1 & 1 & (2A,4D,12B) & 1 & 1 \\
		\hline % inserts single horizontal line
		$5^2:Q_8:3$        & (3A,3B,4A) & 12 & 1 &  &  &  \\
		\hline % inserts single horizontal line
		$5^2:D(2*4)$     & (2A,4A,10C) & 1 & 1 & (2A,4A,10D) & 1 & 1 \\
		& (2B,4A,10A) & 1 & 1 & (2B,4A,10B) & 1 & 1 \\
		& (2A,2B,2C,4A) & 1 & 12 &  &  &  \\
		\hline % inserts single horizontal line
		$5^3:S_5$          & (2A,4A,5G) & 1 & 1 & (2A,4A,5I) & 1 & 1 \\
		& (2A,4A,5H) & 1 & 1 & (2A,4A,5J) & 1 & 1 \\
		\hline % inserts single horizontal line
		$5^3:S_5$          & (2B,4A,6A) & 20 & 1 &  &  &  \\
		\hline % inserts single horizontal line
	\end{tabular}
\end{table}

\begin{table}[H]
	\caption{GOSs for Affine Primitive Groups of Degree 81} % title of Table
	\centering % used for centering table
	\scriptsize\begin{tabular}{|c | c | c| c| c| c| c|} % centered columns (4 columns)
		\hline\hline %inserts double horizontal lines
		group &  ramification type & N.O & L.O &ramification type & N.O & L.O \\
		\hline % inserts single horizontal line
		$3^4:Sp(4,3):2$  & (2C,6N,5A) & 14 & 1&  &  &\\
		\hline % inserts single horizontal line
		$3^4:(GL(1,3) \wr Alt(4))$ &(2D,6H,6N) & 1 & 1& (2D,6G,6M)& 1  & 1\\
		\hline % inserts single horizontal line
		
		$3^4:(Q_8:2)Sym(3)$ &(2A,2C,2D,3F) & 1 & 12& (2A,2C,2D,3G) & 1 &12\\
		\hline % inserts single horizontal line
		$3^4:(2\times Sym(5)) $ &(2C,6H,6N) & 1 & 1&  &  & \\
		\hline % inserts single horizontal line
		$3^4(2\times Alt(6).2)$ &(2C,4C,8B) & 2 & 1& (2C,4C,8D) & 2 &1\\
		&(2C,4B,8A) & 2 & 1& (2C,4B,8C) & 2 &1\\
		\hline % inserts single horizontal line
		$3^4:(Q_8:3):2 $  &(2A,3J,12B) & 1 & 1& (2A,3J,12C) & 1 &1\\
		&(2A,3I,12A) & 1 & 1& (2A,3I,12D) & 1 &1\\
		&(2A,3H,12A) & 1 & 1& (2A,3H,12D) & 1 &1\\
		&(2A,3G,12B) & 1 & 1& (2A,3G,12C) & 1 &1\\
		\hline % inserts single horizontal line
		$3^4:2.Alt(5)$  &(2C,3F,10A) & 1 & 1&(2B,3F,10B)  &1  & 1\\
		&(2B,2C,2C,3F) & 1 & 18&  &  & \\
		\hline % inserts single horizontal line
		$3^4(2^{(3+4)}):4 $ &(2D,4E,8A) & 1 & 1& (2D,4E,8B) & 1 &1\\
		&(2D,4D,8C) & 1 & 1& (2D,4D,8D) & 1 &1\\
		\hline % inserts single horizontal line
		$3^4:Sym(6)$  &(2C,5A,6I) & 6 & 1&  &  &\\
		\hline % inserts single horizontal line
		
		$3^4:2.Alt(5):2$  &(2B,5A,6C)  & 1 & 1&(2B,5A,6B)  & 1 &1\\
		\hline % inserts single horizontal line
		$3^4: (2^3:Alt(4)):Sym(3) $ &(2G,6O,6U) & 3 & 1&(2G,6N,6V)  &3  & 1\\
		
		\hline % inserts single horizontal line
		$3^4: (2^3:2^2):3^2:D_8$   &(2A,6N,8D) & 1 & 1&(2A,6N,8E)  &1  & 1\\
		\hline % inserts single horizontal line
		$3^4: (2^3:2^2):(3^2:4)$  &(2C,4E,8C) & 4 & 1& (2C,4E,8D)  &4  & 1\\
		&(2C,4D,8A) & 4 & 1& (2C,4D,8B)  &4  & 1\\
		
		\hline % inserts single horizontal line
		$3^4: Q_8:Sym(4)$       &(2C,2D,2D,3F) & 1 & 24 &   &  & \\
		
		\hline % inserts single horizontal line
	\end{tabular}
\end{table}

\begin{table}[H]
	\caption{GOSs for Affine Primitive Groups of Degree 27} % title of Table
	\centering % used for centering table
	\scriptsize\begin{tabular}{|c | c | c| c| c| c|c |} % centered columns (4 columns)
		\hline\hline %inserts double horizontal lines
		group &  ramification type & N.O & L.O &ramification type & N.O & L.O \\
		\hline % inserts single horizontal line
		$AGL(3,3)$ &(3E,6D,6F)   & 8  & 1& (3E,6D,6G) & 8 &1\\
		&(3E,4B,4B)   & 16 & 2& (4B,6E,6D) & 8 &1\\
		
		&(4B,4A,6D)   & 12 &  1&  &  & \\
		
		&(2C,6D,26B)  & 1  & 1 & (2C,6D,26C) & 1 & 1\\
		&(2C,6D,26D)  & 1  & 1 & (2C,6D,26A) & 1 & 1\\
		
		&(2B,6F,13B)  & 1 & 1 & (2B,6F,13C) & 1 &1\\
		&(2B,6F,13D)  & 1 & 1 & (2B,6F,13A) & 1 &1\\
		
		&(2B,6G,13B)  & 1 & 1 & (2B,6G,13C) & 1 & 1\\
		&(2B,6G,13D)  & 1 & 1 & (2B,6G,13A) & 1 & 1\\
		
		&(2B,2B,4B,6D) & 1 & 144 & (2B,2B,2C,13A) & 1 &13\\
		&(2B,2B,2C,13B) & 1 & 13 & (2B,2B,2C,13C) & 1 &13\\
		&(2B,2B,2C,13D) & 1 & 13 & (2B,2C,3E,6D)  & 1 & 144\\
		\hline % inserts single horizontal line
		$ASL(3,3)$ &(2A,6B,13A)  & 2 & 1 & (2A,6B,13B) & 2 &1\\
		&(2A,6B,13C)  & 2 & 1 & (2A,6B,13D) & 2 &1\\
		
		&(2A,4A,13A)  & 2 & 1 & (2A,4A,13B) & 2 &1\\
		&(2A,4A,13C)  & 2 & 1 & (2A,4A,13D) & 2 &1\\
		
		&(3F,3F,8A)  & 8  & 2 & (3F,3F,8B) & 8 &2\\
		&(3F,3F,6D)  & 4  & 2 & (3F,3F,6C) & 4 &2\\
		\hline % inserts single horizontal line
		$3^3:S_4\times 2$ &(2A,2E,2E,6H)  & 1  & 12 & (2B,2A,2E,9A)  & 1 & 3\\
		&(2B,2E,2C,6H)  & 1  & 12 & (2B,2B,2E,12A) & 1 & 4\\
		&(2B,2B,2C,2E,2E)  & 1  & 48 & (2B,2A,2E,2E,2E) &1  & 48\\
		&(2A,6A,12A)    &1   & 1  &                &      &\\
		\hline % inserts single horizontal linE
		$3^3:S_4$            &(2A,2A,2B,9A)  & 1  & 3 & (2A,2A,2B,9B) & 1 &3\\
		&(3E,4A,4A)  & 4  & 2 &  &  &  \\
		\hline % inserts single horizontal line
		$3^3:2.A_4$          &(2A,2C,3D,3E)  & 1  & 4 & (2A,6E,9A) & 1  &1 \\
		&(2A,6E,9B)     & 1  & 1 & (2B,6E,6F) & 4  &1 \\
		\hline % inserts single horizontal line
		$3^3:A_4 \times 2$  &(2B,2B,2B,4A)  & 1  & 32 & (2B,4A,9A) & 1 &1\\
		&(2B,4A,9B)  & 1  & 1 &  &  &\\
		\hline % inserts single horizontal line
		$3^3:A(4)$          &(3E,3F,6A)  & 2  & 1 &  &  &\\
		\hline % inserts single horizontal line
	\end{tabular}
\end{table}

\begin{table}[H]
	\caption{GOSs for Affine Primitive Groups of Degree 9} % title of Table
	\centering % used for centering table
	\scriptsize\begin{tabular}{|c | c | c| c| c| c|c |} % centered columns (4 columns)
		\hline\hline %inserts double horizontal lines
		group &  ramification type & N.O & L.O &ramification type & N.O & L.O \\
		\hline % inserts single horizontal line
		$AGL(2,3)$  &(3B,8A,8A)     & 1 & 1  & (3B,8B,8B) & 1 &1\\
		&(3B,6B,8A)     & 1 & 1  & (3B,6B,8B) & 1 &1\\
		
		&(2B,2B,3C,3C)  & 1 & 24 &(2B,2B,3C,4A)& 1 & 36\\
		&(2B,2B,3C,6A)  & 1 & 36 & (2B,2B,4A,6A) & 1 & 48\\
		&(2B,2B,6A,6A)  & 1 & 48 & (2B,2A,3B,8A) & 1 & 4\\
		&(2B,2A,3B,8B)  & 1 & 4  & (2B,3B,3B,6B) & 1 & 12\\
		&(2B,3B,3B,8A)  & 1 & 8  & (2B,3B,3B,8B) & 1 & 8\\
		
		&(2B,2B,3B,3B,3B)  & 1 & 120 & (2B,2B,2A,3B,3B) & 1 & 48\\
		&(2B,2B,2B,2B,3C)  & 1 & 648 & (2B,2B,2B,2B,4A) & 1 & 768\\
		
		&(2B,2B,2B,2B,6A)  & 1 & 864 & (2B,2B,2B,2B,2B,2B) & 1 & 15360\\
		\hline % inserts single horizontal line
		$ASL(2,3)$  &(2A,3C,3C,3D)  & 1 & 3 & (3C,3C,3C,4A) & 1 & 4\\
		&(2A,3B,3C,4A)  & 2 & 2 & (3C,3B,3B,3D) & 1 & 3\\
		&(3B,3C,3C,6B)  & 1 & 2 & (2A,3B,3B,3E) & 1 & 3\\
		
		&(3B,3C,3C,3E)  & 1 & 3 & (3C,3B,3B,6A) & 1 & 2\\
		
		&(3B,3B,3B,4A)  & 1 & 4 & (3E,3D,4A) & 2 & 2\\
		&(3D,3D,6B)     & 1 & 2 & (3E,4A,6B) & 3 & 1 \\
		
		&(6B,6B,6B)     & 2 & 2 & (3E,3E,6A) & 1 & 2 \\
		
		&(3D,4A,6A)     & 3 & 1 & (4A,6B,6A) & 4 & 1 \\
		&(6A,6A,6A)     & 2 & 2 &  &  &  \\
		\hline % inserts single horizontal line
		A$\Gamma$ L(1,9)  &(2B,2B,4B,4B)  & 1 & 32 &  &  &\\
		\hline % inserts single horizontal line
		$ 3^2:D(2*4)$   &(2B,2B,4A,4A)        & 1 & 8  & (2A,2B,3B,4A)       & 1  & 2\\
		&(2A,2B,3A,4A)        & 1 & 2  & (2A,2A,4A,4A)       & 1  & 8\\
		&(2A,2B,2B,2B,4A)     & 1 & 16 & (2A,2A,2B,2B,3B)    & 1  & 4\\
		&(2A,2A,2B,2B,3A)     & 1 & 4  & (2A,2A,2A,2B,4A)    & 1  & 16\\
		&(2A,2A,2B,2B,2B,2B)  & 1 & 32 & (2A,2A,2A,2A,2B,2B) & 1  & 32\\
		\hline % inserts single horizontal line
		AGL(1,9)         &(2A,8B,8A)  & 1 & 1 & (2A,8C,8D) & 1 & 1\\
		\hline % inserts single horizontal line
		$3^2:Q_8=M(9)$    &(4A,4B,4C)  & 4 & 1 &  &  & \\
		\hline % inserts single horizontal line
		$3^2:4$          &(3A,4A,4B)  & 1 & 1 &(3B,4A,4B)  &1  &1 \\
		\hline % inserts single horizontal line
	\end{tabular}
\end{table}

\begin{table}[H]
	\caption{Part1: GOSs for Affine Primitive Groups of Degree 8} % title of Table
	\centering % used for centering table
	\scriptsize\begin{tabular}{|c | c | c| c| c| c|c |} % centered columns (4 columns)
		\hline\hline %inserts double horizontal lines
		group &  ramification type & N.O & L.O &ramification type & N.O & L.O \\
		\hline % inserts single horizontal line
		$AGL(1,8)$  &(2A,7D,7F)  & 1 & 1 & (2A,7B,7C) & 1 &1\\
		&(2A,7A,7E)  & 1 & 1 &            &   &\\
		\hline % inserts single horizontal line
		$A\Gamma L(1,8)$  &(3B,6B,6B)  & 1 & 2 & (3B,6A,7B) & 2 & 1\\
		&(3B,6A,7A)  & 2 & 1 & (3A,7B,7B) & 1 & 2\\
		&(3A,6B,7B)  & 2 & 1 & (3A,6B,7A) & 2 & 1\\
		
		&(2A,3B,3B,3B)  & 1 & 14 & (3A,3A,3B,3B) & 1 & 48\\
		&(2A,3A,3A,3A)  & 1 & 14 &  &  & \\
		\hline % inserts single horizontal line
		$ASL(3,2)$&(3A,7B,7B)  & 2 & 1 &(3A,7A,7B)  &2  &1\\
		&(3A,7A,7A)  & 2 & 1 &(3A,6A,7B)  &2  & 1\\
		&(3A,6A,7A)  & 2 & 1 &(3A,4C,7B)  &2  & 1\\
		&(3A,4C,7A)  & 2 & 1 &(3A,4C,6A)  &4  & 1\\
		&(3A,4A,7A)  & 2 & 1 &(3A,4A,7B)  &2  & 1\\
		&(4B,7A,7B)  & 6 & 1 &(4B,6A,7B)  &4  & 1\\
		&(4B,6A,7A)  & 4 & 1 &(4B,6A,6A)  &2  & 1\\
		&(4B,4C,7A)  & 2 & 1 &(4B,4C,7B)  &2  & 1\\
		
		&(4B,4C,6A)  & 2 & 1 &(4B,4C,4C)  &2  & 2\\
		&(4B,4A,7A)  & 2 & 1 &(4B,4A,7B)  &2  & 1\\
		&(2C,7B,7B)  & 1 & 1 &(2C,7A,7A)  &1  & 1\\
		&(2C,6A,7B)  & 1 & 1 &(2C,6A,7A)  &1  & 1\\
		\hline % inserts single horizontal line
	\end{tabular}
\end{table}

\begin{table}[H]
	\caption{Part2: GOSs for Affine Primitive Groups of Degree 8} % title of Table
	\centering % used for centering table
	\scriptsize\begin{tabular}{|c | c | c| c| c| c|c |} % centered columns (4 columns)
		\hline\hline %inserts double horizontal lines
		group &  ramification type & N.O & L.O &ramification type & N.O & L.O \\
		\hline % inserts single horizontal line
		&(3A,3A,3A,3A)  & 3 & 180 &(3A,3A,3A,4B)  & 2  & 384\\
		&(3A,3A,4B,4B)  & 4 & 288 &(3A,4B,4B,4B)  & 2  & 216\\
		&(4B,4B,4B,4B)  & 3 & 168 &(2C,4B,4B,4B)  & 1  & 72\\
		&(2C,3A,4B,4B)  & 1 & 132 &(2C,3A,3A,4B)  & 1  & 168\\
		
		&(2C,3A,3A,3A)  & 1 & 120 &(2C,2C,3A,4B)  & 1  & 24\\
		&(2C,2C,4B,4B)  & 1 & 24 &(2B,3A,3A,7B)  & 1  & 168\\
		&(2B,3A,3A,7A)  & 1 & 168 &(2B,3A,3A,6A)  & 1  & 240\\
		&(2B,3A,3A,4C)  & 1 & 132 &(2B,3A,3A,4A)  & 1  & 60\\
		
		&(2B,3A,4B,7A)  & 1 & 126 &(2B,3A,3A,7B)  & 1  & 126\\
		&(2B,3A,4B,6A)  & 1 & 168 &(2B,3A,3A,4C)  & 1  & 132\\
		&(2B,3A,4B,4A)  & 1 & 48  &(2B,4B,4B,7B)  & 1  & 105\\
		&(2B,4B,4B,7A)  & 1 & 105 &(2B,4B,4B,6A)  & 1  & 132\\
		
		&(2B,4B,4B,4A)  & 1 & 72 &(2B,4B,4B,4B)  & 1  & 48\\
		&(2B,2C,3A,7A)  & 1 & 42 &(2B,2C,3A,7B)  & 1  & 42\\
		&(2B,2C,3A,6A)  & 1 & 30 &(2B,2C,3A,4C)  & 1  & 24\\
		&(2B,2C,4B,4C)  & 1 & 24 &(2B,2C,4B,7B)  & 1  & 28\\
		&(2B,2C,4B,7A)  & 1 & 28 &(2B,2C,4B,6A)  & 1  & 24\\
		
		&(2B,2C,2C,7A)  & 1 & 7 &(2B,2C,2C,7A)  & 1  & 7\\
		&(2B,2B,7B,7B)  & 2 & 14 &(2B,2B,7A,7A)  & 2  & 14\\
		&(2B,2B,7A,7B)  & 1 & 42 &(2B,2B,6A,7B)  & 1  & 42\\
		
		&(2B,2B,7A,7A)  & 1 & 42 &(2B,2B,6A,6A)  & 1  & 30\\
		&(2B,2B,4C,7B)  & 1 & 28 &(2B,2B,4C,7A)  & 1  & 28\\
		&(2B,2B,4C,6A)  & 1 & 24 &(2B,2B,4C,4C)  & 1  & 24\\
		&(2B,2A,3A,7A)  & 1 & 7 &(2B,2A,3A,7B)  & 1  & 7\\
		&(2B,2A,4B,7A)  & 1 & 7 &(2B,2A,4B,7B)  & 1  & 7\\
		
		&(2A,3A,3A,4B)  & 2 & 14 &(2A,3A,4B,4B)  & 1  & 14\\
		&(2A,4B,4B,4B)  & 2 & 14 &(2B,2B,4A,7B)  & 1  & 14\\
		&(2B,2B,3A,3A,3A)  & 1 & 7812 &(2B,2B,4A,7A)  & 1  & 14\\
		&(2B,2B,3A,3A,4C)  & 1 & 5868 &(2B,2B,3A,4B,4B)  & 1  & 4374\\
		&(2B,2B,4B,4B,4B)  & 1 & 3564 &(2B,2B,2C,3A,3A)  & 1  & 1728\\
		
		&(2B,2B,2C,3A,4B)  & 1 & 1296 &(2B,2B,2C,4B,4B)  & 1  & 912\\
		&(2B,2B,2C,2C,3A)  & 1 & 216 &(2B,2B,2C,2C,4B)  & 1   & 192\\
		&(2B,2B,2B,3A,7A)  & 1 & 1323 &(2B,2B,2B,3A,7B)  & 1  & 1323\\
		&(2B,2B,2B,3A,6A)  & 1 & 1728 &(2B,2B,2B,3A,4A)  & 1  & 1296\\
		&(2B,2B,2B,3A,4B)  & 1 & 432 &(2B,2B,2B,4B,7A)  & 1   & 1029\\
		&(2B,2B,2B,4B,7B)  & 1 & 1029 &(2B,2B,2B,4B,6A)  & 1  & 1296\\
		$ASL(3,2)$       &(2B,2B,2B,4B,4C)  & 1 & 912 &(2B,2B,2B,4B,4A)  & 1   & 384\\
		
		&(2B,2B,2B,2C,7A)  & 1 & 294 &(2B,2B,2B,2C,7B)  & 1  & 294\\
		&(2B,2B,2B,2C,6A)  & 1 & 294 &(2B,2B,2B,2C,4C)  & 1  & 294\\
		&(2B,2B,2B,2A,7A)  & 1 & 49  &(2B,2B,2B,2A,7B)  & 1  & 49\\
		
		&(2B,2B,2A,3A,3A)     & 1 & 210   &(2B,2B,2A,3A,4B)     & 1  & 168\\
		&(2B,2B,2A,4B,4B)     & 1 & 168   &(2B,2B,2B,2B,3A,3A)  & 1  & 60426\\
		
		&(2B,2B,2B,2B,2A,3A)  & 1 & 1512   &(2B,2B,2B,2B,2A,4B)  & 1  & 1344\\
		&(2B,2B,2B,2B,3A,4B)  & 1 & 45360 &(2B,2B,2B,2B,4B,4B)  & 1  & 34992\\
		&(2B,2B,2B,2B,2C,3A)  & 1 & 12960 &(2B,2B,2B,2B,2C,4B)  & 1  & 9600\\
		&(2B,2B,2B,2B,2C,2C)  & 1 & 1680  &(2B,2B,2B,2B,2B,7A)  & 1  & 10290\\
		&(2B,2B,2B,2B,2B,6A)     & 1 & 12960  &(2B,2B,2B,2B,2B,7B)     & 1  & 10290\\
		&(2B,2B,2B,2B,2B,4A)     & 1 & 3360   &(2B,2B,2B,2B,2B,4C)     & 1  & 9600\\
		&(2B,2B,2B,2B,2B,2B,3A)  & 1 & 466560 &(2B,2B,2B,2B,2B,2B,4B)  & 1  & 354240\\
		&(2B,2B,2B,2B,2B,2B,2C)  & 1 & 97920  &(2B,2B,2B,2B,2B,2B,2A)  & 1  & 11760\\
		&(2B,2B,2B,2B,2B,2B,2B,2B) &1 & 3623760   &  &   & \\
		
		\hline % inserts single horizontal line
		
	\end{tabular}
\end{table}

\begin{table}[H]
	\caption{Part1: GOSs for Affine Primitive Groups of Degree 16} % title of Table
	\centering % used for centering table
	\scriptsize\begin{tabular}{|c | c | c| c| c| c|c |} % centered columns (4 columns)
		\hline\hline %inserts double horizontal lines
		group &  ramification type & N.O & L.O &ramification type &N.O & L.O \\
		\hline % inserts single horizontal line
		$2^4:5$                    &(2C,5B,5C)     & 1 & 1  &(2C,5A,5D)     & 1  & 1\\
		&(2B,5B,5C)     & 1 & 1  &(2B,5A,5D)     & 1  & 1\\
		&(2A,5B,5C)     & 1 & 1  &(2A,5A,5D)     & 1  & 1\\
		\hline % inserts single horizontal line
		$2^4:D(2*4)$               &(2D,2D,2C,5B)     & 1 & 5  &(2D,2D,2C,5A)     & 1  & 5\\
		&(2D,2D,2B,5B)     & 1 & 5  &(2D,2D,2B,5A)     & 1  & 5\\
		&(2D,2D,2A,5B)     & 1 & 5  &(2D,2D,2A,5A)     & 1  & 5\\
		&(2D,2D,2D,2D)     & 1 & 60  &(2D,2D,2D,2D,2B)     & 1  & 60\\
		&(2D,2D,2D,2D,2A)     & 1 & 60  &     &   & \\
		
		\hline % inserts single horizontal line
		$2^4:(A(4)\times A(4)):2$  &(3B,6D,6D)     & 1 & 1  &(3B,4A,6C)     & 1  & 1\\
		&(3A,6C,6C)     & 1 & 1  &(3A,4A,6D)     & 1  & 1\\
		&(2C,2C,3B,6B)  & 1 & 2  &(2A,2A,3A,6A)  & 1  & 2\\
		\hline % inserts single horizontal line
		$(2^4:5).4$      &(4C,4C,4A)  & 1 & 1 & (4C,4C,4D) & 2 & 2\\
		&(4B,4C,5A)  & 3 & 1 & (4B,4B,4A) & 1 & 1\\
		&(2C,2C,4B,4C)  & 1 & 18 & (4B,4B,4D) & 2 & 2\\
		\hline % inserts single horizontal line
		$AGL(1,16):2$    &(2B,6B,15A)  & 1 & 1 & (2B,6B,15B) & 1 & 1\\
		&(2B,6A,15C)  & 1 & 1 & (2B,6A,15D) & 1 & 1\\
		\hline % inserts single horizontal line
		$2^4:S(3)\times S(3)$      &(2D,2E,3A,6C)     & 1 & 6  & (2D,2E,3A,4A)   & 1 & 3\\
		&(2C,2E,3A,6B)     & 1 & 6   & (2C,2E,3A,4B)   & 1 & 3\\
		&(2C,2E,3A,4D)     & 1 & 2   & (2C,2E,3A,4C)   & 1 & 1\\
		&(2D,2D,2E,2E,3A)  & 1 & 27  & (2C,2C,2E,2E,3A) & 1 & 27\\
		&(2C,2C,2D,2D,3A)  & 1 & 3   &                  &  &  \\
		\hline % inserts single horizontal line
		$2^4.3^2:4$                &(4C,4C,4D)        & 2 & 2   & (4C,4C,4A)      & 2 & 1\\
		&(4B,4C,6A)        & 3 & 1   & (4B,4B,4D)      & 2 & 2\\
		&(4B,4B,4A)        & 2 & 1   & (3A,4C,8A)      & 1 & 1\\
		&(3A,4B,8B)        & 1 & 1   & (2C,2C,4B,4C)   & 3 & 24\\
		\hline % inserts single horizontal line
		$A\Gamma L(1,16)$            &(4C,4C,6A)        & 3 & 2   & (4B,4B,6A)      & 3 & 2\\
		\hline % inserts single horizontal line
		$S(4)\times S(4):2$        &(4F,6B,6C)        & 6 & 1   & (4F,4E,6B)      & 3 & 1\\
		&(2C,2F,4F,4F)        & 1 & 24   & (2C,2C,4F,8A)      & 1 & 8\\
		&(2C,2C,6C,6C)        & 1 & 12   & (2C,2C,4C,6C)      & 1 & 6\\
		&(2C,2E,4F,6C)        & 1 & 27   & (2C,2E,4F,4E)      & 1 & 12\\
		&(2C,2D,4F,4D)        & 1 & 8    & (2C,2D,4F,4C)      & 1 & 4\\
		&(2C,2D,4F,6A)        & 1 & 6    & (2C,2D,6B,6C)      & 1 & 12\\
		&(2C,2D,4E,6B)        & 1 & 6    & (2C,2D,3A,8A)      & 1 & 2\\
		&(2D,2E,4F,6B)        & 1 & 27   & (2D,2D,4F,4F)      & 1 & 24\\
		&(2D,2D,6B,6B)        & 1 & 12   & (2D,2F,3A,4F)      & 1 & 6\\
		
		&(2C,2C,2D,2F,4F)        & 1 & 48   & (2C,2C,2D,2E,6C)       & 1 & 54\\
		&(2C,2C,2D,2D,4B)        & 1 & 16   & (2C,2C,2D,2E,6D)       & 1 & 24\\
		&(2C,2C,2D,2D,4C)        & 1 & 8   & (2C,2C,2D,2D,6A)        & 1 & 12\\
		&(2C,2D,2E,2E,4F)        & 1 & 120   & (2C,2D,2D,2E,6B)      & 1 & 54\\
		&(2C,2D,2D,2D,4F)        & 1 & 48   & (2C,2D,2D,2F,3A)       & 1 & 12\\
		&(2C,2C,2C,2D,8A)        & 1 & 16   & (2C,2C,2C,2D,2D,2F)      & 1 & 96\\
		&(2C,2C,2D,2D,2E,2E)        & 1 & 240   & (2C,2C,2D,2D,2D,2D)      & 1 & 96\\
		\hline % inserts single horizontal line
		$A\Gamma L(2,4)$             &(4C,4C,6A)        & 3 & 2   & (2C,6C,15A)      & 3 & 1\\
		&(2C,6C,15B)       & 3 & 1   & (3B,6C,6C)       & 3 & 2\\
		&(2C,4C,8A)        & 2 & 1   & (2C,2C,4C,4C)    & 1 & 192\\
		&(2C,2C,3B,5A)     & 1 & 30  & (2C,3B,3B,4C)    & 1 & 32\\
		&(2B,2C,2C,15A)    & 1 & 15  & (2B,2C,2C,15A)   & 1 & 15\\
		&(2B,2C,3B,6C)     & 1 & 42  & (2B,2B,2C,2C,3B) & 1 & 288\\
		\hline % inserts single horizontal line
		$ASL(2,4):2$             &(2C,6A,8A)        & 2 & 1   & (3A,4D,6A)      & 12 & 1\\
		&(4D,4D,5A)        & 1 & 2   & (2B,5A,8A)       & 2 & 1\\
		&(2C,2C,4D,4D)     & 1 & 32  & (2B,2C,2C,8A)    & 2 & 8\\
		&(2B,2C,4D,3A)     & 2 & 30  & (2B,2B,4D,4D)    & 1 & 16\\
		\hline % inserts single horizontal line
	\end{tabular}
\end{table}

\begin{table}[H]
	\caption{Part2: GOSs for Affine Primitive Groups of Degree 16} % title of Table
	\centering % used for centering table
	\scriptsize\begin{tabular}{|c | c | c| c| c| c|c |} % centered columns (4 columns)
		\hline %inserts double horizontal lines
		group &  ramification type & N.O & L.O &ramification type &N.O & L.O \\
		\hline % inserts single horizontal line
		$AGL(2,4)$                   &(2B,6B,15C)        & 1 & 1    & (2B,6B,15B)        & 1 & 1\\
		&(2B,6A,15D)        & 1 & 1    & (2B,6A,15A)        & 1 & 1\\
		&(3D,6C,5B)         & 1 & 1    & (3D,6C,5A)         & 1 & 1\\
		&(3D,3E,15C)        & 1 & 1    & (3D,3E,15B)        & 1 & 1\\
		&(3D,6A,5B)         & 1 & 1    & (3D,6A,5A)         & 1 & 1\\
		&(3C,6D,5B)         & 1 & 1    & (3C,6D,5A)         & 1 & 1\\
		&(3C,3E,15D)        & 1 & 1    & (3C,3E,15A)        & 1 & 1\\
		&(3C,6B,5B)         & 1 & 1    & (3C,6B,5A)         & 1 & 1\\
		&(2B,2B,3D,6C)      & 1 & 18   & (2B,2B,3D,6A)      & 1 & 18\\
		&(2B,2B,3C,6D)      & 1 & 18   & (2B,2B,3D,6B)      & 1 & 18\\
		&(2B,3C,3D,3E)      & 1 & 24   &        &  & \\
		\hline % inserts single horizontal line
		$ASL(2,4)$                   &(3A,3A,5B)         & 6 & 2    & (3A,3A,5A)        & 6 & 2\\
		&(2B,2B,3A,3A)      & 3 & 72   &                   &    & \\
		\hline % inserts single horizontal line
		$2^4.A(6)$                   &(3B,4C,5B)         & 12 & 1     & (3B,4C,5A)        & 12& 1\\
		&(3B,4C,6A)         & 12 & 1     & (4C,4C,5B)        & 4 & 2\\
		&(4C,4C,5A)         & 4  & 2     & (2B,5B,8A)        & 4 & 1\\
		&(3A,5A,5B)         & 2  & 1     & (3A,3B,8A)        & 4 & 1\\
		&(2B,2B,4C,3B)      & 2  & 144   & (2B,2B,4C,4C)     & 3 & 80\\
		&(2B,2B,2B,8A)      & 4  & 24    & (2B,2B,3A,5B)     & 1 & 60\\
		&(2B,2B,3A,5B)      & 1  & 60    & (2B,2B,3A,3B)     & 1 & 36\\
		&(2B,2B,2B,2B,3A)   & 1  & 1728  &                   &   &   \\
		\hline % inserts single horizontal line
		$2^4:S(5)$                   &(5A,6C,6C)         & 3  & 1     & (4B,6C,6C)        & 1 & 1\\
		&(4D,6C,6C)         & 1  & 2     & (6C,6C,6A)        & 1 & 1\\
		&(6C,6C,6B)         & 1  & 2     & (4E,5A,6C)        & 6 & 1\\
		&(4E,4B,6C)         & 1  & 1     & (4E,4D,6C)        & 2 & 1\\
		&(4E,6A,6C)         & 2  & 1     & (4E,6B,6C)        & 4 & 1\\
		&(4E,4E,5A)         & 3  & 1     & (4E,4E,4B)        & 1 & 1\\
		&(4E,4E,6B)         & 4  & 1     & (3A,6C,12A)       & 1 & 1\\
		&(3A,6C,8A)         & 2  & 1     & (3A,4E,12A)       & 2 & 1\\
		&(3A,4E,8A)         & 1  & 1     & (2D,2D,6C,6C)     & 1 & 36\\
		
		&(2D,2D,4E,6C)      & 1  & 54    & (2D,2D,4E,4E)     & 1 & 48\\
		&(2D,2E,3A,6C)      & 1  & 21    & (2D,2E,3A,4E)     & 1 & 30\\
		&(2C,2D,6C,4E)      & 1  & 30    & (2C,2D,6C,4B)     & 1 & 6\\
		&(2C,2D,6C,4D)      & 1  & 12    & (2C,2D,6C,6A)     & 1 & 7\\
		&(2C,2D,6C,6B)      & 1  & 14    & (2C,2D,4E,5A)     & 1 & 45\\
		&(2C,2D,4E,4C)      & 1  & 8     & (2C,2D,4E,4D)     & 1 & 16\\
		&(2C,2D,4E,6A)      & 1  & 10    & (2C,2D,4E,6B)     & 1 & 20\\
		&(2C,2D,3A,12A)     & 1  & 7     & (2C,2D,3A,8A)     & 1 & 18\\
		&(2C,2E,6C,6C)      & 1  & 18    & (2C,2E,4E,6C)     & 1 & 24\\
		&(2C,2E,4E,4E)      & 1  & 24    & (2C,2E,3A,5A)     & 1 & 15\\
		&(2C,2C,5A,5A)      & 1  & 30    & (2C,2C,6C,12A)    & 1 & 6\\
		&(2C,2C,6C,8A)      & 1  & 8     & (2C,2C,4E,12A)    & 1 & 8\\
		&(2C,2C,4E,8A)      & 1  & 8     & (2C,2C,4B,5A)     & 1 & 5\\
		&(2C,2C,4D,5A)      & 1  & 10    & (2C,2C,6A,5A)     & 1 & 5\\
		&(2C,2C,6B,5A)      & 1  & 10    & (2C,2D,2D,2D,6C)  & 1 & 270\\
		
		&(2C,2D,2D,2D,4E)   & 1  & 408   & (2C,2D,2D,2E,3A)  & 1  & 144  \\
		&(2C,2C,2D,2D,5A)   &1   & 225   & (2C,2C,2D,2D,4B)  & 1  & 40\\
		&(2C,2C,2D,2D,4D)   &1   & 80    & (2C,2C,2D,2D,6A)  & 1  & 48\\
		&(2C,2C,2D,2D,6B)   &1   & 96    & (2C,2C,2D,2E,6C)  & 1  & 108\\
		&(2C,2C,2D,2E,4E)   &1   & 144   & (2C,2C,2C,2D,12A) & 1  & 36\\
		&(2C,2C,2C,2D,8A)   &1   & 48    & (2C,2C,2C,2E,5A)  & 1  & 75\\
		&(2C,2C,2C,2D,2D,2E)   &1   & 730   & (2C,2C,2D,2D,2D,2D)  & 1  & 2016\\
		\hline % inserts single horizontal line
		$2^4:A(5)$                   &(3A,5B,5B)         & 3 & 1      & (3A,5A,5B)        & 3 & 1\\
		&(3A,5A,5A)         & 3 & 1      & (3A,5B,6C)        & 1 & 1\\
		&(3A,5B,6A)         & 1  &1      & (3A,5B,6B)        & 1 & 1\\
		&(3A,5A,6B)         & 1  & 1     & (3A,5B,6C)        & 1 & 1\\
		&(3A,5A,6C)         & 1  & 1     & (3A,4A,5A)        & 1 & 1\\
		&(3A,4A,5B)         & 1  & 1     & (3A,4B,5A)        & 1 & 1\\
		&(2C,2C,3A,5B)      & 1  & 45    & (3A,4B,5B)        & 1 & 1\\
		
		&(2C,2C,3A,5A)      & 1  & 45    & (2C,2C,3A,6A)     & 1 & 18\\
		&(2C,2C,3A,6B)      & 1  & 18    & (2C,2C,3A,6C)     & 1 & 18\\
		&(2C,2C,3A,4A)      & 1  & 18    & (2C,2C,3A,4B)     & 1 & 36\\
		&(2C,2C,2C,2C,3A)   & 1  & 810   &                   &   &   \\
		\hline % inserts single horizontal line
	\end{tabular}
\end{table}

\begin{table}[H]
	\caption{Part3: GOSs for Affine Primitive Groups of Degree 16} % title of Table
	\centering % used for centering table
	\scriptsize\begin{tabular}{|c | c | c| c| c| c|c |} % centered columns (4 columns)
		\hline %inserts double horizontal lines
		group &  ramification type & N.O & L.O &ramification type &N.O & L.O \\
		\hline % inserts single horizontal line
		$2^4.A(7)$                   &(3A,6B,6A)         & 6 & 1      & (3A,6B,6B)        & 3 & 2\\
		&(3A,3B,14B)        & 2 & 1      & (3A,3B,14A)       & 2 & 1\\
		&(3A,6B,7B)         & 3  &1      & (3A,6B,7A)        & 3 & 1\\
		&(3A,5A,6A)         & 6  & 1     & (3A,5A,6B)        & 18 & 1\\
		&(3A,5A,7A)         & 8  & 1     & (3A,4A,7B)        & 8 & 1\\
		&(3A,5A,5A)         & 9  & 2     & (3A,4B,14B)       & 6 & 1\\
		&(3A,4B,14A)        & 6  & 1     & (3A,4B,8A)        & 8 & 1\\
		&(3B,4B,7A)         & 6  & 1     & (3B,4B,7B)        & 6 & 1\\
		&(4B,4B,6A)         & 12  & 2    & (4B,4B,6B)        & 18 & 2\\
		&(4B,4B,7B)         & 12  & 2    & (4B,4B,7A)        & 12 & 2\\
		&(4B,4B,5A)         & 32  & 2    & (2B,6B,14B)       & 4 & 1\\
		&(2B,6B,14A)        & 4   & 1    & (2B,7B,14B)       & 2 & 1\\
		&(2B,7B,14A)        & 1   & 1    & (2B,7B,8A)        & 2 & 1\\
		&(2B,7A,14B)        & 1   & 1    & (2B,7A,14A)       & 2 & 1\\
		&(2B,7A,8A)         & 2   & 1    & (2B,5A,14B)       & 4 & 1\\
		&(2B,2B,3A,6A)      & 1  & 132    & (2B,5A,14A)       & 4 & 1\\
		
		&(2B,2B,3A,6B)      & 1  & 270   & (2B,2B,3A,7B)     & 1 & 126\\
		&(2B,2B,3A,7A)      & 1  & 126   & (2B,2B,3A,5A)     & 1 & 450\\
		&(2B,2B,4B,4B)      & 3  & 504   & (2B,2B,2B,14B)    & 2 & 42\\
		&(2B,2B,2B,14A)     & 2  & 42    & (2B,3A,3A,3B)     & 1 & 186\\
		&(2B,3A,3A,4B)      & 1  & 576   & (2B,2B,2B,2B,3A)  & 1  & 9720   \\
		\hline % inserts single horizontal line
		$2^4.S(6)$                   &(4E,4F,6C)         & 8 & 1      & (3B,6B,6C)        & 12 & 1\\
		&(4F ,5A,6B)        &8 & 1       & (4E,6B,6C)        & 6 & 1\\
		&(5A,6B,6B)         & 7  &2      & (2C,5A,12A)        & 2 & 1\\
		&(2E,6C,6C)         & 3  & 2     & (2D,5A,6C)         & 6 & 1\\
		&(2E,6C,12A)        & 6  & 1     & (2B,2B,2B,2B,2E,3B) &1  & 540 \\
		
		&(2C,2C,4F,6B)      & 1  & 12    & (2C,2C,6B,6B)     & 1 & 24\\
		&(2C,2C,3A,5A)      & 1  & 10    & (2C,2D,3A,3B)     & 1 & 12\\
		&(2E,2C,4E,4F)      & 1  & 40    & (2E,2C,3B,6B)     & 2 & 27\\
		&(2E,2C,4E,6B)      & 1  & 42    & (2E,2C,3A,6C)     & 1 & 24\\
		&(2E,2C,2D,5A)      & 1  & 30    & (2E,2E,4F,4F)     & 1  & 32   \\
		&(2E,2E,4F,6B)      & 1  & 120   & (2E,2E,6B,6B)     & 2  & 144   \\
		&(2E,2E,2C,12A)     & 2  & 12    & (2E,2E,2D,6C)     & 1  & 72   \\
		&(2B,2C,3B,5A)      & 2  & 15    & (2B,2C,3B,6A)     & 2  & 6   \\
		&(2B,2C,4F,6C)      & 1  & 18    & (2B,2C,4E,5A)     & 1  & 20   \\
		&(2B,2C,6B,6C)      & 1  & 18    & (2B,2E,3B,6C)     & 2  & 36   \\
		&(2B,2E,4F,5A)      & 1  & 40    & (2B,2E,4E,6C)     & 1  & 42   \\
		&(2B,2E,6B,5A)      & 1  & 70    & (2B,3A,3B,4F)     & 1  & 14   \\
		&(2B,3A,3B,6B)      & 1  & 18    & (2B,2D,3B,4E)     & 1  & 24   \\
		&(2B,2B,6C,6C)      & 1  & 24    & (2B,2B,5A,5A)     & 1  & 20   \\
		&(2B,2B,3B,8A)      & 1  & 8     & (2E,2E,2C,2C,3A)  & 1  & 108   \\
		
		&(2E,2E,2E,2C,2D)   &  1 &   360     & (2B,2E,2C,2C,4F)  & 1  & 80   \\
		&(2B,2E,2E,2C,3B)   &  2 & 171   & (2B,2E,2C,2C,6B)  & 1  & 108   \\
		&(2B,2E,2E,2C,4E)   &  1 & 224   & (2B,2E,2C,2E,4F)  & 1  & 624   \\
		&(2B,2E,2E,2E,6B)   &  1 & 972   & (2B,2B,2C,2C,5A)  & 1  & 50   \\
		&(2B,2B,2C,2D,3B)   &  1 & 54   &  (2B,2B,2E,2C,6C)  & 1  & 108   \\
		&(2B,2B,2E,2E,5A)   &  1 & 350   & (2B,2B,2E,3A,3B)  & 1  & 96   \\
		&(2B,2B,2B,3B,4F)   &  1 & 72    & (2B,2B,2B,3B,6B)  & 1  & 108   \\
		&(2B,2B,2E,2E,2C,2C)   & 1  & 568    & (2B,2B,2E,2E,2E,2E)  & 1  &  5040  \\
		\hline
		$2^4.PSL(4,2)$               &(6B,6B,7B)         & 64  & 1    &   (6B,6B,7A)         &  64 & 1  \\
		&(6B,6B,6C)         & 126  & 1   &   (6B,6B,6A)         &  30 & 1  \\
		&(6B,6B,5A)         & 114  & 1   &   (6B,6B,4E)         &  24 & 1  \\
		
		&(4F,6B,7B)         & 52  & 1    &   (4F,6B,7A)         &  52 & 1  \\
		&(4F,6B,6C)         & 90   & 1   &   (4F,6B,6A)         &  30 & 1  \\
		&(4F,6B,5A)         & 42   & 1   &   (4F,6B,4E)         &  36 & 1  \\
		
		&(4F,4F,7B)         & 18   & 1   &   (4F,4F,7A)         &  18 & 1  \\
		&(4F,4F,4E)         & 24   & 1   &   (4D,6B,7B)         &  18 & 1  \\
		&(4D,6B,7A)         & 18   & 1   &   (4D,6B,6C)         &  24 & 1  \\
		&(4D,6B,5A)         & 36   & 1   &   (4D,4D,6C)         &  4 & 1  \\
		
		&(4D,4D,5A)         & 6    & 1   &   (4C,6C,7B)         &  15 & 1  \\
		&(4C,6C,7A)         & 15   & 1   &   (4C,6C,6C)         &  36 & 1  \\
		
		&(4C,6B,15B)        & 12   & 1   &   (4C,6B,15A)        &  12 & 1  \\
		&(4C,6B,14B)        & 18   & 1   &   (4C,6B,14A)        &  18 & 1  \\
		&(4C,6B,12A)        & 12   & 1   &   (4C,6B,8A)         &  18 & 1  \\
		\hline % inserts single horizontal line
	\end{tabular}
\end{table}

\begin{table}[H]
	\caption{Part4: GOSs for Affine Primitive Groups of Degree 16} % title of Table
	\centering % used for centering table
	\scriptsize\begin{tabular}{|c | c | c| c| c| c|c |} % centered columns (4 columns)
		\hline %inserts double horizontal lines
		group &  ramification type & N.O & L.O &ramification type &N.O & L.O \\
		\hline % inserts single horizontal line
		&(4C,5A,6C)            & 24   & 1      &   (4C,5A,7B)        &  12 & 1  \\
		&(4C,5A,7A)            & 12   & 1      &   (4C,5A,5A)        &  12 & 1  \\
		
		&(4C,6A,6A)            & 6    & 1      &   (4C,6A,5A)        &  12 & 1  \\
		&(4C,4F,15B)           & 12   & 1      &   (4C,4F,15A)       &  12 & 1  \\
		&(4C,4F,14B)           & 14   & 1      &   (4C,4F,14A)       &  14 & 1  \\
		&(4C,4F,12A)           & 18   & 1      &   (4C,4F,8A)        &  12 & 1  \\
		
		&(4C,4E,6C)            & 8    & 1      &   (4C,4E,5A)        &  12 & 1  \\
		&(4C,4D,15B)           & 3    & 1      &   (4C,4D,15A)       &  3 & 1  \\
		
		&(3A,4C,15B)           & 3   & 1       &   (3A,4C,15A)       &  2 & 1  \\
		&(3A,4C,14B)           & 2   & 1       &   (3A,4C,14A)       &  2 & 1  \\
		
		&(3B,6C,7B)            & 11  & 1       &   (3B,6C,7A)        &  11 & 1  \\
		&(3B,6B,15B)           & 14  & 1       &   (3B,6B,15A)       &  14 & 1  \\
		&(3B,6B,14B)           & 16  & 1       &   (3B,6B,14A)        &  16 & 1  \\
		&(3B,6B,12A)           & 10  & 1       &   (3B,6B,8A)        &  10 & 1  \\
		
		&(3B,5A,7B)            & 3  & 1        &   (3B,5A,7A)        &  3 & 1  \\
		
		&(3B,4F,15B)           & 8  & 1        &   (3B,4F,15A)        &  8 & 1  \\
		&(3B,4F,14B)           & 6  & 1        &   (3B,4F,14A)        &  6 & 1  \\
		&(3B,4F,12A)           & 20  & 1       &   (3B,4E,6C)        &  8 & 1  \\
		
		&(3B,4E,5A)            & 4  & 1        &   (3B,4D,15B)        &  4 & 1  \\
		
		&(4F,4D,6C)            & 24 & 1        &   (3B,4D,15A)        &  4 & 1  \\
		&(4F,4D,7A)            & 14  & 1        &  (4F,4D,7B)         &  14 & 1  \\
		&(4F,4D,5A)            & 12  & 1        &                     &   &   \\
		
		&(2D,7B,15B)           & 2  & 1        &   (2D,7A,14B)         &  3 & 1  \\
		&(2D,7B,15A)           & 2  & 1        &   (2D,7B,14A)         &  3 & 1  \\
		&(2D,7B,12A)           & 2  & 1        &   (2D,7A,15B)         &  2 & 1  \\
		&(2D,7A,15A)           & 2  & 1        &   (2D,7A,12A)         &  2 & 1  \\
		&(2D,4E,15B)           & 2  & 1        &   (2D,4E,15A)         &  2 & 1  \\
		
		&(2C,6C,7B)            & 3  & 1        &   (2C,6C,7A)          &  3 & 1  \\
		&(2C,6B,15B)           & 3  & 1        &   (2C,6B,15A)          &  3 & 1  \\
		&(2C,5A,7B)            & 3  & 1        &   (2C,5A,7A)          &  3 & 1  \\
		&(2C,4F,15B)           & 3  & 1        &   (2C,4F,15A)          &  3 & 1  \\
		
		&(2B,15B,15B)          & 1  & 1        &   (2B,15A,15A)         &  1 & 1  \\
		&(2B,14B,15B)          & 1  & 1        &   (2B,14B,15A)         &  1 & 1  \\
		&(2B,14A,15B)          & 1  & 1        &   (2B,14A,15A)         &  1 & 1  \\
		&(2B,12A,15B)          & 1  & 1        &   (2B,12A,15A)         &  1 & 1  \\
		&(2B,8A,15B)           & 1  & 1        &   (2B,8A,15A)         &  1 & 1  \\
		
		&(3A,6B,7A)            & 6  & 1        &  (3A,6B,7B)         &  6 & 1  \\
		&(3A,4D,7A)            & 2  & 1        &  (3A,4D,7B)         &  2 & 1  \\
		
		$2^4.PSL(4,2)$                &(2B,4C,6B,6B)         & 1 & 2700      & (2B,4C,3A,6B)        & 1 & 234\\
		&(2B,4C,4D,6A)         & 1 & 660       & (2B,3A,4C,4D)        & 1 & 48\\
		&(2B,4C,4C,6C)         & 1 & 654       & (2B,4C,4C,5A)        & 1 & 450\\
		&(2B,4C,4F,6B)         & 1 & 2154      & (2B,3A,4C,4F)        & 1 & 252\\
		&(2B,4C,4F,4D)         & 1 & 528       & (2B,4C,4F,4F)        & 1 & 1692\\
		&(2B,2C,6B,6B)         & 1 & 288       & (2B,2C,4C,6C)        & 1 & 90\\
		&(2B,2C,4C,5A)         & 1 & 90        & (2B,2C,4F,6B)        & 1 & 360\\
		&(2B,2C,4F,4F)         & 1 & 288       & (2B,2B,6C,6C)        & 1 & 216\\
		&(2B,2B,6C,7B)         & 1 & 126       & (2B,2B,6C,7A)        & 1 & 126\\
		&(2B,2B,6B,14B)        & 1 & 140       & (2B,2B,6B,14A)       & 1 & 140\\
		&(2B,2B,6B,15B)        & 1 & 120       & (2B,2B,6B,15A)       & 1 & 120\\
		&(2B,2B,6B,12A)        & 1 & 96        & (2B,2B,6B,8A)        & 1 & 120\\
		
		&(2B,2B,5A,7B)        & 1 & 105       & (2B,2B,5A,7A)        & 1 & 105\\
		&(2B,2B,5A,6C)        & 1 & 210       & (2B,2B,5A,5A)        & 1 & 150\\
		&(2B,2B,3A,14B)       & 1 & 14        & (2B,2B,3A,14A)       & 1 & 14\\
		&(2B,2B,3A,15B)       & 1 & 15        & (2B,2B,3A,15A)       & 1 & 15\\
		&(2B,2B,4D,15B)       & 1 & 30        & (2B,2B,4D,15A)       & 1 & 30\\
		&(2B,2B,4E,6C)        & 1 & 60        & (2B,2B,4E,5A)        & 1 & 60\\
		&(2B,2B,6A,6C)        & 1 & 72        & (2B,2B,5A,6A)        & 1 & 60\\
		&(2B,2B,4E,14B)       & 1 & 112       & (2B,2B,4E,14A)       & 1 & 112\\
		&(2B,2B,4E,15B)       & 1 & 90        & (2B,2B,4E,15A)       & 1 & 90\\
		&(2B,2B,4E,12A)       & 1 & 120       & (2B,2B,4E,8A)        & 1 & 96\\
		&(2B,3B,6B,6B)        & 1 & 2448      & (2B,3B,3A,6B)       & 1 & 216\\
		&(2B,3B,4D,6C)       & 1 & 552        & (2B,3B,3A,4D)       & 1 & 60\\
		
		&(2B,3B,4C,6C)       & 1 & 510        & (2B,3B,3A,4C)       & 1 & 450\\
		&(2B,2C,3B,6C)       & 1 & 72         & (2B,2C,3B,3A)       & 1 & 60\\
		\hline
	\end{tabular}
\end{table}

\begin{table}[H]
	\caption{Part5: GOSs for Affine Primitive Groups of Degree 16} % title of Table
	\centering % used for centering table
	\scriptsize\begin{tabular}{|c | c | c| c| c| c|c |} % centered columns (4 columns)
		\hline %inserts double horizontal lines
		group &  ramification type & N.O & L.O &ramification type & N.O & L.O \\
		\hline % inserts single horizontal line
		&(2B,3B,3B,6C)       & 1 & 456        & (2B,3B,3B,3A)       & 1 & 370\\
		&(2B,3B,4F,6B)       & 1 & 1942       & (2B,3B,3A,4F)       & 1 & 198\\
		&(2B,3B,4F,4D)       & 1 & 456        & (2B,3B,4F,4F)       & 1 & 1442\\
		&(2B,2D,6B,6C)       & 1 & 756        & (2B,2D,6B,6A)       & 1 & 198\\
		&(2B,2D,6B,7B)       & 1 & 406        & (2B,2D,6B,7A)       & 1 & 406\\
		&(2B,2D,4E,6B)       & 1 & 168        & (2B,2D,5A,6B)       & 1 & 660\\
		&(2B,2D,3A,7B)       & 1 & 42         & (2B,2D,3A,7A)       & 1 & 42\\
		&(2B,2D,4D,6C)       & 1 & 156        & (2B,2D,4D,5A)       & 1 & 180\\
		&(2B,2D,4D,7B)       & 1 & 98         & (2B,2D,4D,7A)       & 1 & 98\\
		&(2B,2D,4C,14B)      & 1 & 98         & (2B,2D,4C,14A)      & 1 & 98\\
		&(2B,2D,4C,15B)      & 1 & 75         & (2B,2D,4C,15A)      & 1 & 75\\
		&(2B,2D,4C,12A)      & 1 & 84         & (2B,2D,4C,8A)       & 1 & 96\\
		
		&(2B,2D,2C,14B)      & 1 & 15         & (2B,2D,2C,14A)      & 1 & 15\\
		&(2B,2D,3B,14B)      & 1 & 84         & (2B,2D,3B,14A)      & 1 & 84\\
		&(2B,2D,3B,15B)      & 1 & 75         & (2B,2D,3B,15A)      & 1 & 75\\
		&(2B,2D,3B,12A)      & 1 & 66         & (2B,2D,3B,8A)      & 1 & 60\\
		
		&(2B,2D,4F,6C)      & 1 & 540         & (2B,2D,4F,6A)      & 1 & 180\\
		&(2B,2D,4F,7B)      & 1 & 322         & (2B,2D,4F,7A)      & 1 & 322\\
		&(2B,2D,4F,5A)      & 1 & 360         & (2B,2D,4F,4E)      & 1 & 192\\
		&(2B,2D,4C,6C)      & 1 & 252         & (2B,2D,4E,4C)      & 1 & 288\\
		&(2D,4C,4C,6B)      & 1 & 1896        & (2D,3A,4C,4C)      & 1 & 264\\
		&(2D,4C,4C,4D)      & 1 & 424         & (2D,4C,4C,4E)      & 1 & 1680\\
		&(2D,2C,4C,6B)      & 1 & 252         & (2D,2C,4C,4F)      & 1 & 288\\
		&(2D,3B,4C,6B)      & 1 & 1896        & (2D,3B,3A,4C)      & 1 & 138\\
		&(2D,3B,4C,4D)      & 1 & 444         & (2D,3B,4C,4F)      & 1 & 1320\\
		&(2D,2B,2C,6B)      & 1 & 198         & (2D,2B,2C,4F)      & 1 & 180\\
		
		&(2D,3B,3B,6B)      & 1 & 1560        & (2D,3B,3B,4D)      & 1 & 296\\
		&(2D,3B,3B,4F)      & 1 & 648         & (2D,2D,6B,6B)      & 2 & 2088\\
		&(2D,2D,4D,6B)      & 1 & 672         & (2D,2D,4D,4D)      & 1 & 88\\
		&(2D,2D,4C,6B)      & 1 & 618         & (2D,2D,4C,6A)      & 1 & 168\\
		&(2D,2D,4C,7B)      & 1 & 308         & (2D,2D,4C,7A)      & 1 & 308\\
		&(2D,2D,4C,5A)      & 1 & 330         & (2D,2D,4C,4E)      & 1 & 176\\
		&(2D,2D,2C,7B)      & 1 & 42          & (2D,2D,2C,7A)      & 1 & 42\\
		&(2D,2D,3B,7B)      & 1 & 154         & (2D,2D,3B,7A)      & 1 & 154\\
		&(2D,2D,3B,4E)      & 1 & 112         & (2D,2D,4E,6C)      & 1 & 1656\\
		&(2D,2D,4E,4D)      & 1 & 384         & (2B,2D,2D,2D,4C)   & 1 & 10944 \\
		
		$2^4.PSL(4,2)$                &(2B,2B,2B,6B,6B)         & 1 & 23112      & (2B,2B,2B,3A,6B)        & 1 & 1944\\
		&(2B,2B,2B,4D,6B)         & 1 & 5184       & (2B,2B,2B,3A,4D)        & 1 & 432\\
		&(2B,2B,2B,4C,6C)         & 1 & 4698       & (2B,2B,2B,4C,5A)        & 1 & 4050\\
		&(2B,2B,2B,2C,6C)         & 1 & 648        & (2B,2B,2B,2C,5A)        & 1 & 600\\
		&(2B,2B,2B,3B,6C)         & 1 & 4212       & (2B,2B,2B,3B,5A)        & 1 & 3600\\
		&(2B,2B,2B,4F,6B)         & 1 & 18594      & (2B,2B,2B,3A,4F)        & 1 & 1944\\
		&(2B,2B,2B,4F,4D)         & 1 & 4320       & (2B,2B,2B,4F,4F)        & 1 & 14208\\
		&(2B,2B,2B,2D,14B)        & 1 & 784        & (2B,2B,2B,2D,14A)       & 1 & 784\\
		&(2B,2B,2B,2D,12A)        & 1 & 648        & (2B,2B,2B,2D,8A)       & 1 & 672\\
		&(2B,2B,2B,2D,15B)        & 1 & 675        & (2B,2B,2B,2D,15A)       & 1 & 675\\
		&(2B,2B,2D,4C,6B)        & 1 & 16776       & (2B,2B,2D,3A,4C)        & 1 & 1584\\
		&(2B,2B,2D,4C,4D)        & 1 & 3888        & (2B,2B,2D,4C,4F)        & 1 & 13344\\
		
		&(2B,2B,2D,2C,6B)        & 1 & 1944         & (2B,2B,2D,2C,4F)        & 1 & 2016\\
		&(2B,2B,2D,3B,6B)        & 1 & 15084        & (2B,2B,2D,3B,3A)        & 1 & 1140\\
		&(2B,2B,2D,3B,4D)        & 1 & 3264         & (2B,2B,2D,3B,4F)        & 1 & 11784\\
		&(2B,2B,2D,2D,6C)        & 1 & 4536         & (2B,2B,2D,2D,6A)        & 1 & 1296\\
		&(2B,2B,2D,2D,7B)        & 1 & 2548         & (2B,2B,2D,2D,7A)        & 1 & 2548\\
		&(2B,2B,2D,2D,5A)        & 1 & 3900         & (2B,2B,2D,2D,4E)        & 1 & 1152\\
		&(2B,2D,2D,4C,4C)        & 1 & 12120        & (2B,2D,2D,2C,4C)        & 1 & 1728\\
		&(2B,2D,2D,3B,4C)        & 1 & 11544        & (2B,2D,2D,2C,3B)        & 1 & 1296\\
		&(2B,2D,2D,3B,3B)        & 1 & 9804         & (2B,2D,2D,2D,6B)        & 1 & 17064\\
		&(2B,2D,2D,2D,4D)        & 1 & 3840         & (2B,2D,2D,2D,4F)        & 1 & 11232\\
		\multicolumn{7}{|c|}{\textbf{Projection-fiber Algorithm}}\\
		&(2B,2B,2B,2B,2B,6C)         & 1 & 38889       & (2B,2B,2B,2B,2B,5A)        & 1 & 33750\\
		&(2B,2B,2B,2B,2D,6B)         & 1 & 142560      & (2B,2B,2B,2B,2D,3A)        & 1 & $3\times 4320$\\
		&(2B,2B,2B,2B,2D,4D)         & 1 & 31488       & (2B,2B,2B,2B,2D,4F)        & 1 & 113280\\
		&(2B,2B,2B,2D,2D,4C)         & 1 & 106176      & (2B,2B,2B,2D,2D,2C)        & 1 & $3\times 4320$\\
		&(2B,2B,2B,2D,2D,3B)         & 1 & 94824       & (2B,2B,2D,2D,2D,2D)        & 1 & $3\times 34992$\\
		&      &  &       &   (2B,2B,2B,2B,2B,2D,2D)   & 1 & 902400 \\
		\hline
	\end{tabular}
\end{table}

\begin{table}[H]
	\caption{Part1: GOSs for Affine Primitive Groups of Degree 32} % title of Table
	\centering % used for centering table
	\scriptsize\begin{tabular}{|c | c | c| c| c| c|c |} % centered columns (4 columns)
		\hline %inserts double horizontal lines
		group &  ramification type & N.O & L.O &ramification type & N.O & L.O \\
		\hline % inserts single horizontal line
		&(4F,6C,6F)  & 84 & 1 &(4F,4J,12B)  & 72  &1 \\
		&(4F,4J,8B)  & 72 & 1 &(4F,4J,6E)   & 78  &1 \\
		&(4F,4J,5A)  & 78 & 1 &(4F,4I,6F)   & 78  &1 \\
		&(4F,4F,21B) & 10 & 1 &(4F,4F,21A)  & 10  &1 \\
		&(4F,4F,12C) & 12 & 1 &(4F,4F,10A)  & 18  &1 \\
		&(4F,4F,8C)  & 24 & 1 &(3B,4F,14D)  & 10  &1 \\
		&(3B,4F,14C) & 10 & 1 &(3B,4F,6F)   & 48  &1 \\
		
		&(3B,4B,10A) & 6  & 1 &(3B,4B,8C)   & 4  &1 \\
		&(3B,6C,12B) & 84 & 1 &(3B,6C,8B)   & 90  &1 \\
		&(3B,6C,6E)  & 144 & 1 &(3B,6C,6D)   & 12  &1 \\
		&(3B,6C,5A)  & 120 & 1 &(3B,4J,4J)   & 18  &1 \\
		
		&(3B,4I,12B)  & 78  & 1 &(3B,4I,8B)  & 72  &1 \\
		&(3B,4I,6E)   & 90  & 1 &(3B,4I,6D)  & 12  &1 \\
		&(3B,4I,5A)   & 90  & 1 &(3B,4D,8B)  & 12  &1 \\
		&(3B,4D,5A)   & 18  & 1 & (3B,3B,12B)& 48  & 1 \\
		
		$ASL(5,2)$      &(3B,3B,8B)  & 48  & 1 &(3B,3B,6E)   & 48  &1 \\
		&(3B,3B,5A)  & 36  & 1 &(3A,3B,10A)  & 2  &1 \\
		&(3A,3B,8C)  & 2   & 1 &(2E,12B,12B)  & 16  &1 \\
		&(2E,8B,12B)  & 20  & 1 &(2E,8B,8B)  & 16  &1 \\
		&(2E,6E,12B)  & 16  & 1 &(2E,6E,8B)  & 16  &1 \\
		
		&(2E,6E,6E)   & 22  & 1 &(2E,6C,21B)  & 2  &1 \\
		&(2E,6C,21A)  & 2   & 1 &(2E,5A,12B)  & 30  &1 \\
		&(2E,5A,8B)   & 12  & 1 &(2E,5A,6E)   & 18  &1 \\
		&(2E,5A,5A)   & 14  & 1 &(2E,4I,21B)  & 2  &1 \\
		
		&(2E,4F,31A)  & 1   & 1 &(2E,4F,31B)  & 1  &1 \\
		&(2E,4F,31C)  & 1   & 1 &(2E,4F,31D)  & 1  &1 \\
		&(2E,4F,31E)  & 1   & 1 &(2E,4F,31F)  & 1  &1 \\
		
		& (2D,8B,8C )  &  20  & 1 &(2E,3B,15B)  & 1  &1 \\
		&(2E,3B,15A)   & 1   & 1  &(2D,12B,21B)  & 13  &1 \\
		&(2D,12B,21A)  & 13   & 1 &(2D,12B,12C)  & 16  &1 \\
		&(2D,10A,12B)  & 20   & 1 &(2D,8C,12B)   & 20  &1 \\
		
		&(2D,8B,21A)  &  13   & 1 &(2D,8B,21B)  & 13 &1\\
		&(2D,8B,12C)  &  16   & 1 &(2D,8B,10A)  & 20 &1\\
		
		&(2D,6F,14D)   & 13   & 1 &(2D,6F,14C)  & 13  &1 \\
		&(2D,6F,6F)    & 22   & 1 &(2D,6E,21B)  & 18  &1 \\
		&(2D,6E,8C)    & 16   & 1 &(2D,6E,21A)  & 18  &1 \\
		&(2D,6E,12C)   & 16   & 1 &(2D,6E,10A)  & 18  &1 \\
		
		&  (2D,6D,21A)   &  2  & 1 &(2D,6D,21B)  & 2  &1 \\
		&  (2D,5A,21A)   & 18  & 1  &(2D,5A,21B)  & 18  &1 \\
		
		&  (2E,4I,21A) &  2  &  1 &(2D,5A,12C)  & 20  &1 \\
		&(2D,5A,10A)   & 14  & 1  &(2D,5A,8C)   & 12  &1 \\
		&(2D,4J,31A)   & 4   & 1  &(2B,4J,31B)  & 4  &1 \\
		&(2D,4J,31C)   & 4   & 1  &(2B,4J,31D)  & 4  &1 \\
		&(2D,4J,31E)   & 4   & 1  &(2B,4J,31F)  & 4  &1 \\
		
		&(2B,6F,31A)   & 1   & 1  &(2B,6F,31B) & 1  &1 \\
		&(2B,6F,31C)   & 1   & 1  &(2B,6F,31D) & 1  &1 \\
		&(2B,6F,31E)   & 1   & 1  &(2B,6F,31F) & 1  &1 \\
		
		&(2B,10A,21B)  & 1   & 1  &(2B,10A,21A) & 1  &1 \\
		&(2B,8C,21B)   & 1   & 1  &(2B,8C,21A)  & 1  &1 \\
		\hline
	\end{tabular}
\end{table}

\begin{table}[H]
	\caption{Part2: GOSs for Affine Primitive Groups of Degree 32} % title of Table
	\centering % used for centering table
	\scriptsize\begin{tabular}{|c | c | c| c| c| c|c |} % centered columns (4 columns)
		\hline %inserts double horizontal lines
		\multicolumn{7}{|c|}{\textbf{Projection-fiber Algorithm}}\\
		\hline % inserts single horizontal line
		group &  ramification type & N.O & L.O &ramification type & N.O & L.O \\
		\hline % inserts single horizontal line
		
		&(2B,3B,4F,4F) & 1   &432  &(2B,2B,3B,8C) &1  & 24 \\
		&(2B,2E,3B,4I) & 1   &84   &(2B,2B,3B,10A) &1  &30 \\
		&(2D,2E,3B,6C) & 1   & 78  &     &      &    \\
		&(2B,2D,4F,6F) & 1   &588  &(2B,2D,2E,21A) &1  &14 \\
		$ASL(5,2)$       &(2B,2D,2E,12B) & 1   &588 &(2B,2D,2E,21B) &1  &14 \\
		&(2B,2D,2E,6D) & 1   &78   &(2B,2D,2E,8B) &1   &600 \\
		&(2B,2D,2E,6E) & 1   &780  &(2B,2D,2E,5A) &1   &780 \\
		&(2D,2E,3A,3B) & 1   &46   &(2D,2E,3B,4B) &1   &88 \\
		&(2D,2E,4F,4F)  & 1  & 624  &(2D,2D,4F,4J) & 1  & 3360   \\
		&(2D,2D,2E,12B) & 1  &720   &(2D,2D,2E,8B) &1   &672 \\
		&(2D,2D,2E,6E)  & 1  &720   &(2D,2D,2E,5A) &1   &680 \\
		&(2D,2D,2D,21A) & 1  &630   &(2D,2D,2D,21B)&1   &630 \\
		&(2D,2D,2D,12C) & 1  &720   &(2D,2D,2D,8C) &1   &672 \\
		&(2D,2D,2D,10A) & 1  &680   &(2D,2D,3B,4I) &1   &3720 \\
		&(2D,2D,3B,6C)  & 1  &4368  &(2D,2D,3B,4D) &1   & $3\times 88$ \\
		&(2D,2D,3B,3B)  & 1  & 3$\times 560$  &(2B,2B,2D,2E,3B)&1   &528 \\
		&(2B,2D,2D,2D,3B) & 1  & 30024   &(2D,2D,2D,2D,2E) &  1 & 31744 \\
		\hline
	\end{tabular}
\end{table}

\begin{table}[H]
	\caption{GOSs for Affine Primitive Groups of Degree 64 and 128} % title of Table
	\centering % used for centering table
	\scriptsize\begin{tabular}{|c | c | c | c | c | c |c |} % centered columns (4 columns)
		\hline %inserts double horizontal lines
		group &  ramification type & N.O & L.O & ramification type & N.O & L.O \\
		\hline % inserts single horizontal line
		$2^6:3^2:D_{12}$  &(2F,2G,2H,3D) & 1 & 9  &  &   & \\
		\hline % inserts single horizontal line
		$2^6:(3^2:3):4$   &(2D,4F,12A)  & 1 & 1 & (2D,4F,12B) & 1 & 1\\
		&(2D,4E,12C)  & 1 & 1 & (2D,4E,12D) & 1 & 1\\
		\hline % inserts single horizontal line
		$2^6:(3^2:3):D_8$  &(2D,2F,2G,4G)  & 1 & 12 &  &  &\\
		&(2G,4E,12C)    & 1 & 1 & (2G,4E,12D) & 1 & 1\\
		\hline % inserts single horizontal line
		$ 2^6:3^3:D_{12}$   &(2G,6I,6K)        & 3 & 1  &       &   & \\
		\hline % inserts single horizontal line
		$2^6:3^3:Alt(4)$    &(2E,3F,9C)  & 3 & 1 & (2E,3F,9D) & 3 & 1\\
		&(2E,3E,9A)  & 3 & 1 & (2E,3E,9B) & 3 & 1\\
		\hline % inserts single horizontal line
		$2^6:(3^2:3):SD_{16}$    &(2E,4H,8A)  & 3 & 1 &(2E,4H,8B)  &4  &1 \\
		\hline % inserts single horizontal line
		$2^6:3^3:Sym(4)$         &(2F,4G,9A)    & 3 & 1 &(2F,4G,9B)  &3  &1 \\
		&(2F,3D,24A)   & 1 & 1 &(2F,3D,24B)  &1  &1 \\
		&(2F,2F,2F,4G) & 1  & 96  &(2F,2F,2F,3D)  &1  &72 \\
		\hline % inserts single horizontal line
		$2^6:3^2:Sym(4)$    &(2D,4H,9A)  & 1 & 1 &(2D,4H,9B)  &1  &1 \\
		\hline % inserts single horizontal line
		$2^6:(GL(2,2) \wr Sym(3))$  &(2L,4R,12I)  & 4 & 1 &(2E,2L,2L,4S)  & 1 &16 \\
		\hline % inserts single horizontal line
		$2^6:(3^2:3):Q_8:Sym(3)$    &(2E,6C,8C)  & 1 & 1 &(2E,6C,8D)  &1  &1 \\
		\hline % inserts single horizontal line
		$AGL(6,2)$    &(2F,6J,8D)  & 48 & 1 &(2F,6J,6J)  &192  &1 \\
		&(2F,4O,15D)  & 7 & 1 &(2F,4O,15E)  &7  &1 \\
		&(2F,4O,14H)  & 10 & 1 &(2F,4O,14G)  &10  &1 \\
		&(2F,4K,15D)  & 1 & 1 &(2F,4K,15E)  &1  &1 \\
		&(2F,4H,21A)  & 1 & 1 &(2F,4H,21B)  &1  &1 \\
		&(2F,3C,42A)  & 2 & 1 &(2F,3C,42B)  &2  &1 \\
		&(2F,3C,21A)  & 4 & 1 &(2F,3C,21B)  &4  &1 \\
		&(2D,8D,7E)  & 16 & 1 &(2D,6J,7E)  &82  &1 \\
		&(2D,6F,14H)  & 4 & 1 &(2D,6F,14G)  &4  &1 \\
		&(2B,7E,15D)  & 1 & 1 &(2B,7E,15E)  &1  &1 \\
		&(3C,4O,6F)  &  92 & 1 &  &  & \\
		\hline % inserts single horizontal line
		$A\Sigma L(3,4)$    &(2C,4D,14A)  & 1 & 1 &(2C,4D,14B)  &1  &1 \\
		\hline % inserts single horizontal line
		$2^6:3.Alt(6)$ &(3C,3D,4D)        & 6 & 1 &            &     & \\
		\hline % inserts single horizontal line
		$2^6:(3\times GL(3,2))$    &(3C,3E,4D)  & 4 & 1 &  &  & \\
		\hline % inserts single horizontal line
		$2^6:Sp(6, 2)$    &(2H,6H,7A)  & 14 & 1 &(2H,6I,7A)  &42  &1 \\
		\hline % inserts single horizontal line
		$2^6:GO-(6, 2)$    &(2I,6G,8F)  & 4 & 1 &(2I,4R,12G)   &4  &1 \\
		&(2I,4R,10B)  & 5 & 1 &(2I,4R,8E)   &6  &1 \\
		&(2I,4P,10B)  & 3 & 1 &(2I,4O,12I)  &4  &1 \\
		&(2I,4O,9A)  & 6 & 1 &(2I,4L,12I)   &2  &1 \\
		&(2I,4L,9A)  & 3 & 1 &(2C,8F,12I)   &1  &1 \\
		&(2C,8F,9A)  & 1 & 1 &(2I,2I,2I,10B)  & 1 &192 \\
		&(2I,2I,2I,4L)  & 1 & 96  &(2F,2I,2I,4R)  & 1 & 112\\
		&(2F,2I,2I,4P)  & 1 & 60 &(2C,2I,2I,8F)  & 1 & 32 \\
		\hline % inserts single horizontal line
		$2^6:O-(6,2)$    &(2E,4H,12F) & 3 & 1 &(2E,4H,12E)  & 3  &1 \\
		&(2E,4H,9A)  & 3 & 1 &(2E,4H,9B)  & 3  &1 \\
		&(2E,4F,9A)  & 3 & 1 &(2E,4F,9B)  & 3  &1 \\
		\hline % inserts single horizontal line
		$2^6:Sym(8)$     &(2J,6G,7A)   & 6 & 1 &(2J,4S,12F)  & 4  &1 \\
		&(2J,4S,10B)  & 6 & 1 &(2J,4R,8E)  & 6  &1 \\
		&(2J,4O,15A)  & 4 & 1 &(2J,2J,2J,4O)  & 1  &192 \\
		&(2F,2J,2J,4S)  & 1 & 48 &(2D,2J,2J,4R)  & 1  & 88\\
		\hline % inserts single horizontal line
		$2^6:Alt(8)$     &(2G,4G,15A)  & 1 & 1 &(2G,4G,15B)  & 1  &1 \\
		&(2G,4H,15A)  & 1 & 1 &(2G,4H,15B)  & 1  &1 \\
		&(2C,4J,15A)  & 1 & 1 &(2C,4J,15B)  & 1  &1 \\
		\hline % inserts single horizontal line
		$2^6:Sym(7)$     &(2H,4N,12F)  & 2 & 1 &(2H,4N,10B)  & 3  &1 \\
		&(2H,4M,12G)  & 4 & 1 &(2H,3B,12I)  & 2  &1 \\
		&(2H,3B,12H)  & 1 & 1 &(2D,7A,12B)  & 1  &1 \\
		&(2H,2H,2I,3B)  &1  & 66 &(2H,2H,2H,4L)  &  1 & 48\\
		&(2F,2H,2H,4N)  & 1 &  48&  &   & \\
		\hline % inserts single horizontal line
		$2^6:Alt(7)$     &(2E,4E,7A)  & 6 & 1 &(2E,4E,7B)  & 6  &1 \\
		\hline % inserts single horizontal line
		$2^6:\Sigma U(3,3)$     &(2D,6A,6C)  & 6 & 1 &  &   & \\
		\hline % inserts single horizontal line
		$2^6:SU(3,3)$      &(2B,6A,7A)  & 2 & 1 &(2B,6A,7B)  & 2  &1 \\
		\hline % inserts single horizontal line
		$2^6:PGL(2,7)$     &(2G,4I,8A)  & 2 & 1 &(2G,4I,8B)  & 2  &1 \\
		&(2G,4G,8A)  & 1 & 1 &(2G,4G,8B)  & 1  &1 \\
		&(2G,4G,12A) & 1 & 1 &(2G,2G,2F,3A)  & 1  & 30 \\
		&(2D,2G,2G,4I)  & 1 & 24 &(2D,2G,2G,4G)  & 1  & 12\\
		\hline % inserts single horizontal line
		$ASL(7,2)$         &(2F,3B,14K)  & 1 & 1 &(2F,3B,14L)  & 1  & 1\\
		\hline % inserts single horizontal line
	\end{tabular}
\end{table}

%\section{GAP-Codes}
%\begin{singlespace}
%	\include{append}
%\end{singlespace}
\begin{verbatim}
GAP-Codes
LoadPackage ("mapclass");
# Generate ramification types that fit the Riemann Hurwitz number
# Assign the group g first, then fix the list of conjugacy class representatives
CC:=[];
ct:=[];
CheckingTheGroup:=function(group)
CC:=List(ConjugacyClasses(group),Representative);
ct:=CharacterTable(group);
end;

# The index of a permutation
PermIndex:=function(perm,degree)
return degree-Length(Orbits(Group(perm),[1..degree]));
end;

# The indices of conjugacy class representatives for a group
FindInd:=function(group,degree)
local i,IndexSet,t,Ind;
IndexSet:=[];
for i in [2..Length(CC)] do
t:=CC[i];
Ind:=PermIndex(t,degree);
Append(IndexSet,[rec(pos:=i,index:=Ind , Ord:=Order(t),
size:=Size(Centralizer(k,t)))]);
od;
return IndexSet;
end;

# Find all possible ramification types for a group with fixed degree and genus
Dim:=[];
RamiTypes:=function(group,degree,genus)
local h,bb,dim,n,RH,IndexSet,Indices,PossibleCombinations,
i,a,b,Temp,j,k,c,RamificationTypes;
dim:=Length(Factors(degree));
n:=Filtered(NormalSubgroups(group),x->Size(x)=degree)[1];
RamificationTypes:=[];
RH:=2*(degree+genus-1);
IndexSet:=FindInd(group,degree);
Indices:=List(IndexSet,x->x.index);
PossibleCombinations:=RestrictedPartitions(RH,Elements(Indices)); 
for i in [1..Length(PossibleCombinations)] do
Temp:=[];
a:=PossibleCombinations[i];
b:=Elements(a);
for j in [1..Length(b)] do
k:=Length(Filtered(a,x->x=b[j]));
c:=List(Filtered(IndexSet,x->x.index=b[j]),x->x.pos);
Append(Temp,[UnorderedTuples(c,k)]);
od;
Temp:=Cartesian(Temp);
for j in [1..Length(Temp)] do
Append(RamificationTypes,[Concatenation(Temp[j])]);
od;
od;
for i in [1..Length(RamificationTypes)] do
a:=RamificationTypes[i];
b:=List(a,x->Size(Centralizer(n,CC[x])));
c:=[];
for j in [1..Length(b)] do
if b[j]<> 1 then
Append(c,[dim-Length(FactorsInt(b[j]))]) ;
else
Append(c,[dim]);
fi;
od;
bb:=List(a,x->CC[x]);
h:=List(bb,x->degree-Length(MovedPoints(x)));
if Sum(c)<2*dim then
Unbind(RamificationTypes[i]);
elif Sum(c)=2*dim and not 0 in h then
Unbind(RamificationTypes[i]);
else
Append(Dim,[c]);
fi;
od;
return Elements(RamificationTypes); 
end;

AddOneGenerator:=function(SubgroupList,NewGen,group)
local NewSubgroupList,a,b,x,Cgx,gg,c,t,y,j,r,h,flag,hh;
NewSubgroupList:=[];
for a in [1..Length(CC)] do
NewSubgroupList[a]:=[];
od; 
for b in [1..Length(CC)] do
x:=CC[b];
Cgx:=Centralizer(group,x);
for gg in SubgroupList[b] do
for c in DoubleCosetRepsAndSizes(group,Centralizer(group,NewGen),
Normalizer(Cgx,gg)) do
t:=NewGen^(c[1]);
y:=x*t;
j:=1;
while j < Length(CC)+1 do
r:=RepresentativeAction(group,y,CC[j]);
if not r=fail then 
break; 
else 
j:=j+1;
fi;
od;
h:=Group(Concatenation(GeneratorsOfGroup(gg),[t]))^r;
flag:=0; 
for hh in NewSubgroupList[j] do
if IsConjugate(Centralizer(group,CC[j]),hh,h) then
flag:=1; 
break;
fi;
od;
if flag=0 then
Add(NewSubgroupList[j],h);
fi;
od;
od;
od;
return NewSubgroupList;
end;

AddOneGenerator1:=function(SubgroupList,NewGen,group,p)
local NewSubgroupList,a,b,x,Cgx,gg,c,t,y,r,h,stop;
stop:=0;
for b in [1..Length(CC)] do
x:=CC[b];
Cgx:=Centralizer(group,x);
for gg in SubgroupList[b] do
for c in DoubleCosetRepsAndSizes(group,Centralizer(group,NewGen),
Normalizer(Cgx,gg)) do
t:=NewGen^(c[1]);
y:=x*t;
r:=RepresentativeAction(group,y,CC[p]);
if r<>fail then 
h:=Group(Concatenation(GeneratorsOfGroup(gg),[t]))^r;
if h=group then
stop:=1;
break;
fi;
fi;
od;
od;
od;
return stop;
end;

# Find generating types
GeneratingType:=function(group,degree,genus)
local RamificationTypes,GeneratingTypes,i,SubgroupList,
ClassRepTuple,NewGen,m,j,k,n,a,p;
GeneratingTypes:=[];
RamificationTypes:=RamiTypes(group,degree,genus);
for i in [1..Length(RamificationTypes)] do
Print("\r","Checking the ramification type ",i," with ",
Length(RamificationTypes)-i," remaining ","\c");
SubgroupList:=[];
p:=0;
ClassRepTuple:=List(RamificationTypes[i],x->CC[x]);
for k in [1..Length(CC)] do
SubgroupList[k]:=[];
if IsConjugate(group,CC[k]^-1,
ClassRepTuple[Length(ClassRepTuple)]) then
p:=k;
fi;
od;

SubgroupList[RamificationTypes[i][1]]:=[Group(ClassRepTuple[1])];
for j in [2..Length(ClassRepTuple)-1] do 
NewGen:=ClassRepTuple[j];
if j=Length(ClassRepTuple)-1 then
m:=AddOneGenerator1(SubgroupList,NewGen,group,p);
if m=1 then
Append(GeneratingTypes,[RamificationTypes[i]]);
fi;
else
m:=AddOneGenerator(SubgroupList,NewGen,group);
fi;
SubgroupList:=m;
od;
od;
n:=Concatenation("CasesFor","degree",String(degree),"group",
String(Position(AllPrimitiveGroups(DegreeOperation,degree),k)),
"genus",String(genus));
if Length(GeneratingTypes)<>0 then
AppendTo(n,"GeneratingTypes:=",GeneratingTypes,";\n");
AppendTo(n,"group:=",group,";\n");
AppendTo(n,"CC:=",CC,";\n");
fi;
Print("\n");
return GeneratingTypes;
end; 

# Lifting Generating braid orbits
LiftingQuotientorbit:=function(group,degree,tuple)
local GeneratingTypes,kk,Qorbits,gg,g,genus,q,f,l,n,e,nn,phi,s,ss,OO,orbits,
xx,i,j,A,LL,PP,a,b,c,d,h,m,z,zz,qz;

Qorbits:=[];
kk:=Stabilizer(group,degree);;
gg:=GeneratorsOfGroup(group);;
q:=Socle(group);;
nn:=List(gg,x->x*RepresentativeAction(q,degree^x,degree));;
phi:=GroupHomomorphismByImages(group,kk,gg,nn);;
s:=tuple;
ss:=List(s,x->x^phi);
orbits:=GeneratingMCOrbits(kk,0,ss);;
for j in [1..Length(orbits)] do
xx:=orbits[j].TupleTable[1].tuple;;
LL:=[];
for i in [1..Length(xx)] do
A:=Elements(PreImages(phi,xx[i]));;
A:=Filtered(A,x->IsConjugate(group,x,s[i]));
Add(LL,A);
od;
PP:=[];
if Length(LL)=3 then
for a in LL[1] do
for b in LL[2] do
c:=(a*b)^-1;
if c in LL[3] and Subgroup(group,[a,b])=group then 
Add(PP,[a,b,c]);
fi;
od;
od;
if Size(Center(kk))=1 then
OO:=Orbits(q,PP,function(p,e) return [p[1]^e,p[2]^e,p[3]^e]; end);;
l:=Length(OO);
else
zz:=Center(kk).1;
z:=Random(PreImages(phi,zz));
qz:=Subgroup(group,Concatenation(GeneratorsOfGroup(q),[z]));
OO:=Orbits(qz,PP,function(p,e) return [p[1]^e,p[2]^e,p[3]^e]; end);;
l:=Length(OO);
fi;
fi;
if Length(LL)=4 then
for a in LL[1] do
for b in LL[2] do
for c in LL[3] do
d:=(a*b*c)^-1;
if d in LL[4] and Subgroup(group,[a,b,c])=group then 
Add(PP,[a,b,c,d]);
fi;
od;
od;
od;
if Size(Center(kk))=1 then
OO:=Orbits(q,PP,function(p,e) return [p[1]^e,p[2]^e,p[3]^e,p[4]^e]; end);;
l:=Length(OO);
else
zz:=Center(kk).1;
z:=Random(PreImages(phi,zz));
qz:=Subgroup(group,Concatenation(GeneratorsOfGroup(q),[z]));
OO:=Orbits(qz,PP,function(p,e) return [p[1]^e,p[2]^e,p[3]^e,p[4]^e]; end);;
l:=Length(OO);
fi;
fi;
if Length(LL)=5 then
for a in LL[1] do
for b in LL[2] do
for c in LL[3] do
for d in LL[4] do
e:=(a*b*c*d)^-1;
if e in LL[5] and Subgroup(group,[a,b,c,d])=group then 
Add(PP,[a,b,c,d,e]);
fi;
od;
od;
od;
od;
if Size(Center(kk))=1  then
OO:=Orbits(q,PP,function(p,e) return [p[1]^e,p[2]^e,p[3]^e,p[4]^e,p[5]^e];end);;
l:=Length(OO);
else
zz:=Center(kk).1;
z:=Random(PreImages(phi,zz));
qz:=Subgroup(group,Concatenation(GeneratorsOfGroup(q),[z]));
OO:=Orbits(qz,PP,function(p,e) return [p[1]^e,p[2]^e,p[3]^e,p[4]^e,p[5]^e];end);;
l:=Length(OO);
fi;
fi;
if Length(LL)=6 then
for a in LL[1] do
for b in LL[2] do
for c in LL[3] do
for d in LL[4] do
for e in LL[5] do
f:=(a*b*c*d*e)^-1;
if f in LL[6] and Subgroup(group,[a,b,c,d,e])=group then 
Add(PP,[a,b,c,d,e,f]);
fi;
od;
od;
od;
od;
od;
if Size(Center(kk))=1  then
OO:=Orbits(q,PP,function(p,e) return [p[1]^e,p[2]^e,p[3]^e,p[4]^e,p[5]^e,
p[6]^e];end);;
l:=Length(OO);
else
zz:=Center(kk).1;
z:=Random(PreImages(phi,zz));
qz:=Subgroup(group,Concatenation(GeneratorsOfGroup(q),[z]));
OO:=Orbits(qz,PP,function(p,e) return [p[1]^e,p[2]^e,p[3]^e,p[4]^e,p[5]^e,
p[6]^e];end);;
l:=Length(OO);
fi;
fi;
if Length(LL)=7 then
for a in LL[1] do
for b in LL[2] do
for c in LL[3] do
for d in LL[4] do
for e in LL[5] do
for f in LL[6] do
g:=(a*b*c*d*e*f)^-1;
if g in LL[7] and Subgroup(group,[a,b,c,d,e,f])=group then 
Add(PP,[a,b,c,d,e,f,g]);
fi;
od;
od;
od;
od;
od;
od;
if Size(Center(kk))=1 then
OO:=Orbits(q,PP,function(p,e) return [p[1]^e,p[2]^e,p[3]^e,p[4]^e,p[5]^e
,p[6]^e,p[7]^e];end);;
l:=Length(OO);
else
zz:=Center(kk).1;
z:=Random(PreImages(phi,zz));
qz:=Subgroup(group,Concatenation(GeneratorsOfGroup(q),[z]));
OO:=Orbits(qz,PP,function(p,e) return [p[1]^e,p[2]^e,p[3]^e,p[4]^e,p[5]^e
,p[6]^e,p[7]^e];end);;
l:=Length(OO);
fi;
fi;
if Length(LL)=8 then
for a in LL[1] do
for b in LL[2] do
for c in LL[3] do
for d in LL[4] do
for e in LL[5] do
for f in LL[6] do
for g in LL[7] do
h:=(a*b*c*d*e*f*g)^-1;
if h in LL[8] and Subgroup(group,[a,b,c,d,e,f,g])=group then 
Add(PP,[a,b,c,d,e,f,g,h]);
fi;
od;
od;
od;
od;
od;
od;
od;
if Size(Center(kk))=1 then
OO:=Orbits(q,PP,function(p,e) return [p[1]^e,p[2]^e,p[3]^e,p[4]^e,p[5]^e
,p[6]^e,p[7]^e,p[8]^e];end);;
l:=Length(OO);
else
zz:=Center(kk).1;
z:=Random(PreImages(phi,zz));
qz:=Subgroup(group,Concatenation(GeneratorsOfGroup(q),[z]));
OO:=Orbits(qz,PP,function(p,e) return [p[1]^e,p[2]^e,p[3]^e,p[4]^e,p[5]^e
,p[6]^e,p[7]^e,p[8]^e];end);;
l:=Length(OO);
fi;
fi;
if Length(LL)=9 then
for a in LL[1] do
for b in LL[2] do
for c in LL[3] do
for d in LL[4] do
for e in LL[5] do
for f in LL[6] do
for g in LL[7] do
for h in LL[8] do
m:=(a*b*c*d*e*f*g*h)^-1;
if m in LL[9] and Subgroup(group,[a,b,c,d,e,f,g,h])=group then 
Add(PP,[a,b,c,d,e,f,g,h,m]);
fi;
od;
od;
od;
od;
od;
od;
od;
od;
if Size(Center(kk))=1 then
OO:=Orbits(q,PP,function(p,e) return [p[1]^e,p[2]^e,p[3]^e,p[4]^e,p[5]^e
,p[6]^e,p[7]^e,p[8]^e,p[9]^e];end);;
l:=Length(OO);
else
zz:=Center(kk).1;
z:=Random(PreImages(phi,zz));
qz:=Subgroup(group,Concatenation(GeneratorsOfGroup(q),[z]));
OO:=Orbits(qz,PP,function(p,e) return [p[1]^e,p[2]^e,p[3]^e,p[4]^e,p[5]^e
,p[6]^e,p[7]^e,p[8]^e,p[9]^e];end);;
l:=Length(OO);
fi;
fi;
if l<>1 then
Add(Qorbits,[rec(numberofquotientorbit:=j, 
LargestLength:=Length(orbits[j].TupleTable), q:=l, List:=PP )]);
else
Add(Qorbits,[rec(numberofquotientorbit:=j, 
LargestLength:=Length(orbits[j].TupleTable), q:=l)]);
fi;
od;
n:=Concatenation("CasesFor","degree",String(degree),"group");
if Length(Qorbits)<>0 then
AppendTo(n,"group:=",group,";\n");
AppendTo(n,"GT:=",GT,";\n"); 
AppendTo(n,"Lifting orbits:=",Qorbits,";\n");
fi;
Print("\n");
return Qorbits;
end; 

# Checking ramification type
Find3Tuple:=function(RamificationType,group)
local GeneratingTuples,ClassRepTuple,g1,g2,g3,gg,Cgx,c,t,y,r,h,i,TT,flag;
ClassRepTuple:=List(RamificationType,x->CC[x]);
GeneratingTuples:=[];
g1:=ClassRepTuple[1];
g2:=ClassRepTuple[2];
g3:=ClassRepTuple[3];
gg:=Group(g1);
Cgx:=Centralizer(group,g1);
for c in DoubleCosetRepsAndSizes(group,Centralizer(group,g2),Cgx) do
t:=g2^(c[1]);
y:=g1*t;
r:=RepresentativeAction(group,y,Inverse(g3));
if r<>fail then                                 
h:=Group(Concatenation(GeneratorsOfGroup(gg),[t]));
if h=group then
TT:=[g1,t,Inverse(y)];
flag:=0;
for i in [1..Length(GeneratingTuples)] do
if RepresentativeAction(group,TT,GeneratingTuples[i],OnTuples)<>fail then
flag:=1;
break;
fi;
od;
if flag=0 then
Append(GeneratingTuples,[TT]);
fi;
fi;
fi;
od;
return GeneratingTuples;
end;

# Checking ramification type in quotient group
QFind3Tuple:=function(tuple,group,degree)
local GeneratingTuples,ClassRepTuple,g1,g2,g3,gg,
Cgx,c,t,y,r,h,i,TT,flag,Cl,gp,kk,q,nn,phi;
kk:=Stabilizer(group,degree);;
gp:=GeneratorsOfGroup(group);;
q:=Socle(group);;
nn:=List(gp,x->x*RepresentativeAction(q,degree^x,degree));;
phi:=GroupHomomorphismByImages(group,kk,gp,nn);;
Cl:=tuple;
ClassRepTuple:=List(Cl,x->x^phi);
GeneratingTuples:=[];
g1:=ClassRepTuple[1];
g2:=ClassRepTuple[2];
g3:=ClassRepTuple[3];
gg:=Group(g1);
Cgx:=Centralizer(kk,g1);
for c in DoubleCosetRepsAndSizes(kk,Centralizer(kk,g2), Cgx) do
t:=g2^(c[1]);
y:=g1*t;
r:=RepresentativeAction(kk,y,Inverse(g3));
if r<>fail then                                 
h:=Group(Concatenation(GeneratorsOfGroup(gg),[t]));
if h=kk then
TT:=[g1,t,Inverse(y)];
flag:=0;
for i in [1..Length(GeneratingTuples)] do
if RepresentativeAction(kk,TT,GeneratingTuples[i],OnTuples)<>fail then
flag:=1;
break;
fi;
od;
if flag=0 then
Append(GeneratingTuples,[TT]);
fi;
fi;
fi;
od;
return GeneratingTuples;
end;

# Find affine primitive groups of degree p^e.
AffinePrimitiveGroups:=function(degree)
local 
GroupL,GroupLL,k,i;
GroupL:=[];
GroupLL:=AllPrimitiveGroups(DegreeOperation,degree);;
for i in [1..Length(GroupLL)] do
k:=GroupLL[i];
if Size(Socle(k))=degree  then
Append(GroupL, [k]);
fi;
od;
return GroupL;
end;

# Find affine primitive group of degree p^e.
AffinePrimitiveGroup:=function(degree,pos)
local 
GroupL,position;
GroupL:=[];
position:=AffinePrimitiveGroups(degree);
if pos>Length(position) then
Print("There are only ", Length(position) ," affine primitive groups");
else
GroupL:=position[pos];
fi;
return GroupL;
end;

# Ordering the conjugacy class representives.
Ordering:=function(group,degree)
local i,j,IndexSet,x,y,oo,c,d,a,Label,LL,zz;
IndexSet:=[];
Label:=FindInd(group,degree);
for i in [1..Length(Label)] do
for j in [i..Length(Label)] do
if i<>j then
if Label[i].Ord> Label[j].Ord then
zz:= Label[i];
Label[i]:= Label[j];
Label[j]:=zz;
fi;
if Label[i].Ord = Label[j].Ord  and Label[i].size < Label[j].size then
zz:= Label[i];
Label[i]:= Label[j];
Label[j]:=zz;
fi;
if Label[i].Ord = Label[j].Ord  and Label[i].size = Label[j].size 
and Label[i].index > Label[j].index then
zz:= Label[i];
Label[i]:= Label[j];
Label[j]:=zz;
fi;
if Label[i].Ord = Label[j].Ord  and Label[i].size = Label[j].size 
and Label[i].index = Label[j].index then
oo:=Set(IndexSet,x->x.Ord);
x:=Label[i].con;
y:=Label[j].con;
for c in [2..Length(CC)-1] do
for d in [c+1..Length(CC)] do
for a in oo do
if IsConjugate(group,x^a,CC[c])=true and 
IsConjugate(group,y^a,CC[d])=true then
if Size(Centralizer(group,CC[c])) < 
Size(Centralizer(group,CC[d])) then 
zz:= Label[i];
Label[i]:= Label[j];
Label[j]:=zz;
else
break;
fi;
if Size(Centralizer(group,CC[c]))=Size(Centralizer(group,CC[d]))
and PermIndex(CC[c],degree) > PermIndex(CC[d],degree) then 
zz:= Label[i];
Label[i]:= Label[j];
Label[j]:=zz;
else
break;
fi;
else
break;
fi;
od;
od;
od;
fi;
fi;
od;
Add(IndexSet, Label[i]);
od;
LL:=[];
for i in [1..Length(IndexSet)] do
x:=IndexSet[i].pos;
y:=IndexSet[i].Ord;
Append(LL,[rec(pos:=x,order:=y)]);
od;
return LL;
end;

# Computing the conjugacy class by random conjugation 
to check a triple (aa,bb,cc).
h:=Centralizer(group,aa);;
Size(h);
Index(group,Centralizer(group,bb));
orbs:=[];
orbs:=List(orbs,x->[x,Order(aa*x)]);;
for i in [1..1000] do
x:=bb^Random(group);
o:=Order(aa*x);;
new:=true;
for i in [1..Length(orbs)] do
if orbs[i][2]=o then
if RepresentativeAction(h,orbs[i][1],x)<>fail then
new:=false;
break;
fi;
fi;
od;
if new then
Add(orbs,[x,o]);
fi;
od;
sum:=0;
for i in [1..Length(orbs)] do
sum:=sum+Index(h,Centralizer(h,orbs[i][1]));
od;
sum;
# We stop if 
sum=Index(group,Centralizer(group,bb))
goodorbs:=Filtered(orbs,x->RepresentativeAction(group,(aa*x[1])^-1,cc)<>fail);;
Length(goodorbs);
GT:=Filtered(goodorbs,x->Size(group)=Size(Group(aa,x[1])));;
Length(GT);

# Application of Lemma 3.13 in [20].
Minusidentity:=function(group,degree,genus,p,e)
local m,q,I,C,n,ll,U,i,j,B,P,PP,RamificationTypes,A; 
n:=Filtered(NormalSubgroups(group),x->Size(x)=degree)[1];
RamificationTypes:=RamiTypes(group,degree,genus);
if p=3 and e=3 then
m:=[[1,0,0],[0,1,0],[0,0,1]]*Z(p);;
q:=Permutation(m,AsList(GF(p)^e));
fi;
if p=3 and e=4 then
m:=[[1,0,0,0],[0,1,0,0],[0,0,1,0],[0,0,0,1]]*Z(p);;
q:=Permutation(m,AsList(GF(p)^e));
fi;
if p=3 and e=5 then
m:=[[1,0,0,0,0],[0,1,0,0,0],[0,0,1,0,0],[0,0,0,1,0],[0,0,0,0,1]]*Z(p);;
q:=Permutation(m,AsList(GF(p)^e));
fi;
if p=3 and e=6 then 
m:=[[1,0,0,0,0,0],[0,1,0,0,0,0],[0,0,1,0,0,0],
[0,0,0,1,0,0],[0,0,0,0,1,0],[0,0,0,0,0,1]]*Z(p);
q:=Permutation(m,AsList(GF(p)^e));
fi;
if p=5 and e=3 then
m:=[[1,0,0],[0,1,0],[0,0,1]]*Z(p)^2;;
q:=Permutation(m,AsList(GF(p)^e));
fi;
if p=5 and e=4 then
m:=[[1,0,0,0],[0,1,0,0],[0,0,1,0],[0,0,0,1]]*Z(p)^2;;
q:=Permutation(m,AsList(GF(p)^e));
fi;
if p=7 and e=3 then
m:=[[1,0,0],[0,1,0],[0,0,1]]*Z(p)^3;;
q:=Permutation(m,AsList(GF(p)^e));
fi;
if p=7 and e=4 then
m:=[[1,0,0,0],[0,1,0,0],[0,0,1,0],[0,0,0,1]]*Z(p)^3;;
q:=Permutation(m,AsList(GF(p)^e));
fi;
P:=[];
for i in [1..Length(RamificationTypes)] do
U:=List(RamificationTypes[i],x->CC[x]);;
ll:=[U[1]*q,U[2]*q,U[3]];;
B:=List(ll,x->Size(Centralizer(n,x)));
A:=[];
for j in [1..Length(B)] do
if B[j]<> 1 then
Append(A,[e-Length(Factors(B[j]))]);
else
Append(A,[e]);
fi;
od;
if Sum(A)>=2*e then
Add(P,RamificationTypes[i]);
fi;
od;
PP:=[];
for i in [1..Length(P)] do
U:=List(P[i],x->CC[x]);;
ll:=[U[1]*q,U[2],U[3]*q];;
B:=List(ll,x->Size(Centralizer(n,x)));
A:=[];
for j in [1..Length(B)] do
if B[j]<> 1 then
Append(A,[e-Length(Factors(B[j]))]);
else
Append(A,[e]);
fi;
od;
if Sum(A)>=2*e then
Add(PP,P[i]);
fi;
od;
C:=[];
for i in [1..Length(PP)] do
U:=List(PP[i],x->CC[x]);;
ll:=[U[1],U[2]*q,U[3]*q];;
B:=List(ll,x->Size(Centralizer(n,x)));
A:=[];
for j in [1..Length(B)] do
if B[j]<> 1 then
Append(A,[e-Length(Factors(B[j]))]);
else
Append(A,[e]);
fi;od;
if Sum(A)>=2*e then
Add(C,PP[i]);
fi;
od;
return C;
end;
\end{verbatim}

%\renewcommand{\bibname}{List of References}
%\addcontentsline{toc}{chapter}{\quad\,\,{List of References}}
%\begin{singlespace}
	%\bibliographystyle{plain}
	% Include the bibliography file without its .bib extension
%	\bibliography{has}
%\end{singlespace}

\end{document}